\documentclass[11pt]{article}
\usepackage{latexsym}
\usepackage{amsfonts,amssymb,amsmath,mathtools,mathabx}
\usepackage{bbold}           
\usepackage{authblk}

\newtheorem{thm}{Theorem}[section]
\newtheorem{cor}[thm]{Corollary}
\newtheorem{lem}[thm]{Lemma}
\newtheorem{pro}[thm]{Proposition}

\newcommand{\ov }{\overline }
\newcommand{\minus}{\smallsetminus}
\newcommand{\x}{\hspace{-0.025in}\times\hspace{-0.025in}}
\newcommand{\e}{\varepsilon}

\title{On the complexity of the word problem \\    
of the R.\! Thompson group $V$}     

\author{J.C.\ Birget}

\date{\scriptsize
13 v 2026}           

\begin{document} 
\maketitle 

\vspace{-0,5cm}

\hspace{4,0cm} Dedicated to Mikhail V.\ Volkov on his 70th birthday

\begin{abstract}
We analyze Lehnert and Schweitzer's proof that the word problem of the
Thompson group $V$ is co-context-free and show that this word problem is 
the complement of the cyclic closure of a union of reverse deterministic 
context-free languages. 
The same is true for any finitely generated subgroup of $V$. 
For certain finite generating sets of $V,$ the word problem is the 
complement of the cyclic closure of the union of four deterministic 
context-free languages. Therefore the word problem of $V$ has 
{\em quadratic} time-complexity on a deterministic multitape Turing machine, 
and belongs to {\sf logDCFL}.
\end{abstract}

{\footnotesize Keywords: Word problem, R.\ Thompson's group $V,$ complexity, 
deterministic context-free languages.

 }

\section{Introduction}

We will use the following notation; see e.g.\ \cite{HU, Ginsburg, 
Harrison, Sipser}.  For an alphabet $A$, the set of all words over $A$ is 
denoted by $A^*$; this includes the empty word $\e$. We let  $A^+$ $=$ 
$A^* \minus \{\e\}$.
The length of $w \in A^*$ is denoted by $|w|$. A {\em language over} $A$ is 
any subset of $A^*$. Rigorously, a language is a pair $(A, L)$ where $A$ is
a finite alphabet and $L \subseteq A^*$; then the complement of a language 
$L$ over $A$ is unambiguously defined as $A^* \minus L$. A word in
$\{0,1\}^*$ is also called a {\em bitstring}. We only use finite alphabets.  

\noindent The following sets of languages are used:

\medskip

${\sf DTime}(T)$, or ${\sf DTime}(T(n))$, is the set of languages 
accepted by deterministic multitape Turing 

\hspace{1,0cm}  machines with {\em time-complexity} function $\,\le T(.)$
 \ \cite{HU,Sipser};

\smallskip

{\sf CFL} is the set of {\em context-free} languages 
\cite{HU, Harrison, Ginsburg, Sipser};

\smallskip

{\sf coCFL} is the set of {\em co-context-free} languages, i.e., 
the languages with context-free complement; 

\smallskip

{\sf DCFL} is the set of {\em deterministic context-free} languages
\cite{HU, Harrison, GinsGrei66};

\smallskip

${\sf DCFL}^{\rm rev}$ is the set of {\em reverse deterministic 
context-free} languages, i.e., the languages whose 

\hspace{1,0cm} reverse is in {\sf DCFL};

\smallskip

$\cup_{\ell\,} {\sf DCFL}$ is the set of languages that are the union of
$\,\le \ell\,$ {\sf DCFL} languages \cite{GinsGrei66, Yamakami};

\smallskip

$\cup {\sf DCFL}$ is the set of languages that are a union of finitely
many {\sf DCFL} languages.

\smallskip

\noindent See our {\sf Appendix} for some details about (deterministic) 
push-down automata, {\sf CFL}, {\sf coCFL}, and {\sf DCFL}. 

All the finite alphabets that we use are subsets of some fixed countable
set; therefore the {\em class} of all finitely generated groups, and all 
complexity classes such as ${\sf DTime}(T)$, {\sf CFL}, {\sf coCFL},
{\sf DCFL}, etc., are {\em sets}. 

This paper is an updated version of \cite{JCBV}.

\subsection{Overview}

{\bf Background:}
The group $V$ of Richard Thompson is a well known finitely presented 
infinite simple group \cite{Th0,Th,CFP,Hig74,BiThomps}. 
Lehnert and Schweitzer \cite{LehSch} proved that the word problem of the 
Higman-Thompson groups $G_{n,r}$, and in particular the group $V$ 
($= G_{2,1}$), over any finite monoid generating set, is co-context-free. 
Hence the word problem of $V$ is in $\,{\sf coNTime}(n)$. And it is in 
${\sf DTime}(n^{2.38})$, 
using Valiant's algorithm for {\sf CFL} and fast boolean matrix 
multiplication; there exist slightly smaller, and more complicated, upper 
bounds than $2.38\,$ (see the literature on the complexity of matrix 
multiplication; the best bound so far is still $> 2.37$). 
Valiant's algorithm can be implemented with the above time-complexity on a 
multitape Turing machine \cite[{\small Summary}]{Valiant}.
As ${\sf DTime}(T(n))$ is closed under complementation, this time-complexity
also holds for {\sf coCFL}.
It had previously been proved that the word problem of $V$ is in log-space 
uniform ${\sf AC}^1$ and in ${\sf DTime}(n^3)\,$ \cite{BiThomps}. 
By the results of Muller and Schupp \cite{MullSchupp83}, context-free word 
problems of groups are in ${\sf DTime}(n)$, since they are actually in 
{\sf DCFL}. Before \cite{LehSch}, the role of {\sf coCFL} for the word 
problem had been studied in \cite{HRRT}.
The classes {\sf CFL} and {\sf coCFL} are subclasses of {\sf logCFL} (the 
class of languages reducible to {\sf CFL} languages by many-one log-space 
reductions); {\sf logCFL} is closed under complementation, and has a nice 
circuit characterisation that explicitly places {\sf logCFL} between 
log-space uniform ${\sf NC}^1$ and ${\sf AC}^1$; see \cite{Johnson}.

\bigskip

\noindent {\bf Results:}
 \ We prove that the complement of the word problem of $V,$ over any finite 
monoid generating set $A$, is the cyclic closure of a union of reverse 
deterministic context-free languages; i.e., 

\smallskip

\hspace{2,0cm} ${\sf wp}_A(V)$  $\,\in\,$ 
${\sf co}({\sf cyc}(\cup_{\ell\,} {\sf DCFL}^{\rm rev}))$
 \ \ $(\,\subseteq {\sf coCFL})\,$,

\smallskip

\noindent where $\,\ell$ depends on the chosen generating set $A$.
 
The same holds for the word problem ${\sf wp}_B(\langle B \rangle)$ of any 
subgroup $\langle B \rangle \subseteq V$, with finite generating set $B\,$ 
(with $\ell$ now depending on $B$).

For a certain finite generating set $\Gamma_{\!H}$ of $V$ (namely Higman's
generating set of four involutions), 

\smallskip

\hspace{2,0cm} ${\sf wp}_{\Gamma_{\!H}}(V)$ $\,=\,$ 
$({\sf wp}_{\Gamma_{\!H}}(V))^{\rm rev}$
$\,\in\,$ $\,{\sf co}({\sf cyc}(\cup_{4\,}{\sf DCFL}))$
  \ \ $(\,\subseteq {\sf coCFL})\,$.

\smallskip

More generally, if a subgroup of $V$ is generated by a finite set $B$ 
of {\em involutions} then $(.)^{\rm rev}$ can be dropped, and we have (for 
some $\,\ell$ depending on $B$):
 
\smallskip

\hspace{2,0cm} ${\sf wp}_B(\langle B \rangle) \,=\,$
$({\sf wp}_B(\langle B \rangle))^{\rm rev}$  $\,\in\,$
$\,{\sf co}({\sf cyc}(\cup_{\ell\,}{\sf DCFL}))$
  \ \ $(\,\subseteq {\sf coCFL})\,$.

\smallskip

Since $\,\cup{\sf DCFL} \subseteq {\sf DTime}(n)$, which is closed under 
reversal and complementation, and since cyclic closure increases 
time-complexity by a factor $n$, it follows that the word problem of $V$ 
over any finite generating set is in $\,{\sf DTime}(n^2)$. The same holds
for all finitely generated subgroups of $V$.

We prove some closure properties of {\sf logDCFL} and 
${\sf log}_{(1:1)}{\sf DCFL}$, that imply that the word problems of $V$
and its finitely generated subgroups are in ${\sf log}_{(1:1)}{\sf DCFL}$
$(\,\subseteq {\sf logDCFL})$. Here ${\sf log}_{(1:1)}$ denotes one-one 
log-space reduction. 

Since all the Higman-Thompson groups $G_{n,r}$ are subgroups of $V,$ the 
results also apply to $G_{n,r}$. 

\smallskip

We give a short introduction to the Thompson group at the end of Section 2.

\subsection{More definitions}

For an alphabet $A$, let $A^{-1}$ be a (not necessarily disjoint)
copy of $A$. The elements of $A^{-1}$ are called the {\em inverse letters}.
Inversion of letters is treated notationally as an involution, i.e., 
$(a^{-1})^{-1}$ denotes $a$. We denote $A \cup A^{-1}$ by $A^{\pm 1}$.
A {\em group generating set} of a group $G$ is a subset $A \subseteq G$ 
such that every element of $G$ can be expressed as the product of a sequence 
of elements of $A^{\pm 1}$.
A {\em monoid generating set} of a group $G$ is a subset $A \subseteq G$ such
that every element of $G$ can be expressed as the product of a sequence of
elements of $A$. Hence, if $A$ is a group generating set then $A^{\pm 1}$ is
a monoid generating set.  We only use {\em internal} generating sets; i.e., 
$A = A^{\pm 1} \subseteq G$, and $a^{-1}$ is the inverse of $a$ in $G$. 

Let $G$ be a finitely generated group, with finite monoid generating set $A$.
If $u, v \in A^*$ represent the same element of $G$ we denote this by 
$\,u =_G v$. The {\bf word problem} of $G$ over $A$ is defined by

\smallskip 

\hspace{1,0cm} ${\sf wp}_{A\!}(G) \ = \ \{w \in A^* :\, w =_G \e \}$.

\smallskip

\noindent The word problem is a language over the alphabet $A$, but the 
languages that arise as word problems of groups have special properties. 
E.g., ${\sf wp}_{A\!}(G)$ is closed under Kleene star, under conjugation, 
and in particular under cyclic permutation; cyclic closure also plays a 
role in conjugacy (see e.g., \cite{CHHR}).

\section{Cyclic closure, complement, reversal, union}

\subsection{Definitions and basic properties of languages}

\noindent The {\em cyclic closure} of a word $w \in A^*$, or of a language 
$L \subseteq A^*$, or of a set ${\cal C}$ of languages, is defined by 

\medskip

\hspace{1,0cm} ${\sf cyc}(w) \,=\, \{ yx \,:\,$
$x, y \in A^* \ {\rm and} \ w = xy\}$, 

\medskip

\hspace{1,0cm} ${\sf cyc}(L) \,= \,\bigcup_{w \in L} {\sf cyc}(w)$, 
 \ \ and 
 \ \ ${\sf cyc}({\cal C}) \,=\,\{{\sf cyc}(L) : L \in {\cal C}\}$.

\medskip 

\noindent This is also called the closure under cyclic permutation. 

A language $L$ is called {\em cyclically closed} \ iff \ $L = {\sf cyc}(L)$.  
A set $\cal C$ of languages is called cyclically closed \ iff 
 \ ${\sf cyc}({\cal C}) \,\subseteq\, {\cal C}$.

\bigskip

\noindent For any set $\cal C$ of languages we define

\smallskip

\hspace{1,0cm} ${\sf co}{\cal \,C}$ $\,=\,$  $\{A_{_L}^* \minus L \,:\,$
$L \in {\cal C}, \, A_{_L}$ is the alphabet of $L\}$.

\smallskip

\noindent In particular, the {\em co-word-problem} of a group $G$ over the
finite monoid generating set $A$ is defined by

\medskip

\hspace{1,0cm} ${\sf cowp}_{A\!}(G)$  $\,=\,$
$A^* \minus {\sf wp}_{A\!}(G)$.

\medskip

\noindent (About the spelling: We write ``co-word-problem'', and not
``co-word problem'', because we have no such thing as a ``co-word''.)

\bigskip

\noindent The {\em reverse} $\,w^{\rm rev}$ of a word $w \in A^*$ is 
defined by induction on length as follows: 

\smallskip

\hspace{1,0cm} $(v a)^{\rm rev} = a\,v^{\rm rev}$ \ for all $a \in A$ and 
$v \in A^*$; \ \ and \ \ $\e^{\rm rev} = \e$.

\smallskip

\noindent  For a language $L \subseteq A^*$, or a set ${\cal C}$ of 
languages, we define 

\smallskip

\hspace{1,0cm} $L^{\rm rev} \,=\, \{w^{\rm rev} :$ $w \in L\}$, \ \ and 
 \ \ ${\cal C}^{\rm rev} \,=\, \{L^{\rm rev} : L \in {\cal C}\}$.

\smallskip

\noindent A language $L$ is called {\em closed under reversal} \ iff
 \ $L = L^{\rm rev}$. A set $\cal C$ of languages is called closed under 
reversal \ iff \ ${\cal C}^{\rm rev} \,\subseteq\, {\cal C}$.
 \ (About terminology: {\em reversal} is an action, the result of which is 
the {\em reverse}; i.e., applying reversal to $L$ produces the reverse of 
$L$.  Compare with inversion versus inverse.)

\begin{lem}\!\!\!. \label{LEMcycClos} 

\smallskip

\noindent {\small \sf (1.a)} If $L \subseteq A^*$ is cyclically closed 
then so is the complement $\,A^* \minus L$. \ I.e., for every 
$L \subseteq A^*:$ 

\smallskip

\hspace{1,0cm}
 ${\sf cyc}({\sf co}({\sf cyc}(L))) \,=\, {\sf co}({\sf cyc}(L))$.

\smallskip

\noindent {\small \sf (1.b)} For all $L \subseteq A^*:$
 \ \ $A^* \minus {\sf cyc}(L) \ \subseteq \ A^* \minus L$ 
$ \ \subseteq \ $ ${\sf cyc}(A^* \minus L)$.

 \ \ In other words, \ ${\sf co}({\sf cyc}(L)) \ \subseteq \ {\sf co}(L)$
$ \ \subseteq \ {\sf cyc}({\sf co}(L))$.
 
 \ \ But usually, $\,{\sf cyc}(A^* \minus L)$  $\,\ne\,$ 
$A^* \minus {\sf cyc}(L)$.
 \ I.e., $\,{\sf co}(.)$ and $\,{\sf cyc}(.)$ do {\em not} commute.

\smallskip

\noindent {\small \sf (2)} \ For all $L_1, L_2 \subseteq A^*:$ 
 \ ${\sf cyc}(L_1 \cup L_2) \,=\, {\sf cyc}(L_1) \cup {\sf cyc}(L_2)$.  
 \ So $\,{\sf cyc}(.)$ and $\,\cup(.)$ commute. 

\smallskip

\noindent {\small \sf (3)} \ For all $L, L_1, L_2 \subseteq A^*:$
 \ \ $\,(A^* \minus L)^{\rm rev} = A^* \minus L^{\rm rev}$, \ and
 \ \ $(L_1 \cup L_2)^{\rm rev} = L_1^{\rm rev} \cup L_2^{\rm rev}$.
\end{lem}
{\sc Proof.} (1.a) Here, $L = {\sf cyc}(L)$.  Let $A^* \minus L = \ov{L}$.  
If $L = A^*$, then $\ov{L} = \varnothing$ is closed under ${\sf cyc}(.)$. 
If $L \ne A^*$, consider $w \in \ov{L}$.  
If ${\sf cyc}(w) \not\subseteq \ov{L}$ then there exists 
$u \in {\sf cyc}(w) \cap L$. Since $L$ is closed under ${\sf cyc}(.)$, this 
implies ${\sf cyc}(u) \subseteq L$, hence 
${\sf cyc}(w) = {\sf cyc}(u) \subseteq L$, hence $w \in L$. 
This contradicts $w \in \ov{L}$. 

\noindent (1.b) For any $L \subseteq A^*$, $L \subseteq {\sf cyc}(L)$ 
implies $A^* \minus {\sf cyc}(L) \subseteq A^* \minus L$. And 
$A^* \minus L \subseteq {\sf cyc}(A^* \minus L)$.

\noindent For example for $A = \{a,b\}$ and $L = \{ab\}$ we have 
$\,\{a,b\}^* \minus {\sf cyc}(ab) = \{a,b\}^* \minus \{ab, ba\}$
$\,\ne\,$  $\{a,b\}^* = {\sf cyc}(\{a,b\}^* \minus \{ab\})$. 

\noindent (2) If $x \in {\sf cyc}(L_1 \cup L_2)$ then there exists 
$u \in L_1 \cup L_2$ such that $x \in {\sf cyc}(u)$. If $u \in L_1$ then
$x \in {\sf cyc}(L_1)$; if $u \in L_2$ then $x \in {\sf cyc}(L_2)$. So,
$x \in {\sf cyc}(L_1) \cup {\sf cyc}(L_2)$.
The converse is straightforward since $L_1 \subseteq L_1 \cup L_2$ implies 
${\sf cyc}(L_1) \subseteq {\sf cyc}(L_1 \cup L_2)$, and similarly for $L_2$.
 
\noindent (3) is straightforward.  
 \ \ \  \ \ \ $\Box$

\bigskip

\noindent The following is straightforward and well known:
For every group $G$ with monoid generating set $A$,

\smallskip

\hspace{1,0cm} ${\sf wp}_{A\!}(G)\,$ is cyclically closed.

\smallskip

\noindent From this and Lemma \ref{LEMcycClos}(1.a) we obtain:
For every group $G$ with monoid generating set $A$,

\smallskip

\hspace{1,0cm} ${\sf cowp}_{A\!}(G)\,$ is cyclically closed.

\begin{lem} \label{LEMcycrevComm} {\bf ($\,{\sf cyc}(.)$ and $(.)^{\rm rev}$
 commute).}

\medskip

For all $L \subseteq A^*:$ \ \ \ \ ${\sf cyc}(L^{\rm rev})$ $\,=\,$
$({\sf cyc}(L))^{\rm rev}$.
\end{lem}
{\sc Proof.} For any $x \in A^*$: $\,x \in {\sf cyc}(L^{\rm rev})$ \ iff 
 \ there exists $u \in {\sf cyc}(x)$ such that $u \in L^{\rm rev}$. 
This means that for some $\alpha,$ $\beta$ $\in$ $A^*$: 
$\,x = \beta \alpha$ and $u = \alpha \beta \in L^{\rm rev}$. 
Equivalently, $x^{\rm rev} = \alpha^{\rm rev} \beta^{\rm rev}$ and 
$u^{\rm rev} = \beta^{\rm rev} \alpha^{\rm rev} \in L$, which is equivalent
to $\,u^{\rm rev} = v \in {\sf cyc}(x^{\rm rev})$ such that $v \in L$.
This means $x^{\rm rev} \in {\sf cyc}(L)$, i.e., 
$\,x \in ({\sf cyc}(L))^{\rm rev}$.
  \ \ $\Box$

\bigskip

\noindent {\bf Some more properties of languages:}

\smallskip

It is a non-trivial fact that the set {\sf CFL} {\em is closed under} 
${\sf cyc}(.)$, i.e., $\,{\sf cyc}({\sf CFL}) \subseteq {\sf CFL}\,$ 
(\cite{Oshiba}, \cite{Maslov}, and \cite[{\small solved Ex.\! 6.4c}]{HU});
hence ${\sf co}({\sf cyc}({\sf CFL})) \subseteq {\sf coCFL}$.
This fact plays an essential role in \cite{LehSch}.

\smallskip

It is easy to prove that {\sf CFL} is closed under $(.)^{\rm rev}\,$ (by 
using grammars); hence, {\sf coCFL} is closed under $(.)^{\rm rev}$.
The class {\sf DCFL} is not closed under $(.)^{\rm rev}\,$ and not closed
under union \cite{HU, Harrison}. 

\smallskip

For all $\ell \ge 1$, $ \ \cup_{\ell\,} {\sf DCFL}$ $\,\subsetneqq\,$ 
$\cup_{\ell+1} {\sf DCFL}$ \ (this is the {\sf DCFL} {\em union hierarchy});
and $\,\bigcup_{\ell = 1}^{\infty} \cup_{\ell\,} {\sf DCFL}$
$=$ $\cup {\sf DCFL}$ $\,\subsetneqq\,$ ${\sf CFL}$ \ (see e.g.\ 
\cite{Yamakami}). 
The class {\sf DCFL} is closed under complementation \cite{GinsGrei66}; 
but $\cup {\sf DCFL}$ is not.  

\smallskip

{\sf DCFL} is {\sl not} closed under ${\sf cyc}(.)$. An example is 
$L = \{a c^n d c^n : n \ge 1\} \,\cup\, \{b c^n d c^{2n} : n \ge 1\}$;
then $\,{\sf cyc}(L) \cap \{c,d\}^* \!\cdot\! \{a, b\} = $
$\{c^n d c^n a : n \ge 1\} \,\cup\, \{c^n d c^{2n} b : n \ge 1\}$, which is
not in {\sf DCFL}; hence, ${\sf cyc}(L)$ is not in {\sf DCFL}.  Here we use 
the fact that the intersection of any $L \in {\sf DCFL}$ with a finite-state 
language is in ${\sf DCFL}\,$ \cite{GinsGrei66}. 
 \ More generally the following is known: 

\smallskip

\begin{pro}\!\!\!. \label{DCFLcyc}
 \ \ {\sf DCFL}, $\,\cup_{\ell\,} {\sf DCFL}$ and $\cup {\sf DCFL}$, are 
{\em not} closed under cyclic permutation.
\end{pro}
{\sc Proof.} E.g., $L_0 = \{w c w^{\rm rev} : w \in \{a,b\}^*\}$
$\in$ {\sf DCFL}. But ${\sf cyc}(L_0) \cap c\{a,b\}^*$ $=$ 
$c\,\{w w^{\rm rev} : w \in \{a,b\}^*\}$ $\,\not\in\,$ $\cup {\sf DCFL}$,
for the same reason as $\,\{w w^{\rm rev} : w \in \{a,b\}^*\}$ $\,\not\in\,$
$\cup {\sf DCFL}$. The latter was stated by Ginsburg and Greibach 
\cite{GinsGrei66}, and proved by Yamakami 
\cite[{\small Thm.\ 1.5(1)}]{Yamakami}.
 \ \ \  \ \ \ $\Box$

\medskip

\subsection{The word problem and reversal}

\noindent The following result might be known, but it is hard to find any 
mention in the literature.

\begin{pro} \label{PROPnotrev} {\bf (non-closure of the word problem 
under reversal).} 
 
\smallskip

\noindent For any finitely generated group $G$ the following are 
equivalent:

\smallskip

$\bullet$ For {\em every} finite monoid generating set $A$ of $G$, 
the set $\,{\sf wp}_{A\!}(G)$ is closed under reversal. 

\smallskip

$\bullet$ The group $G$ is commutative.

\medskip

\noindent Hence every non-commutative finitely generated group has a 
finite monoid generating set for which the word problem is {\em not} 
closed under reversal. 
\end{pro}
{\sc Proof.} If $G$ is commutative, then ${\sf wp}_{A\!}(G)$ is 
obviously closed under reversal.

Conversely, let $G$ be any group with finite generating set $A$. If $G$ is
1-generated, it is commutative. So let us take the case where $|A| \ge 2$. 
For any two generators $a, b \in A$ we have $ab b^{-1} a^{-1} =_G \e$. 
By a Tietze transformation we can add a new generator $c$ and the relation
$c = b^{-1} a^{-1}$, i.e., $abc =_G \e$.
 \ If $\,{\sf wp}_{A \cup \{c\}\!}(G)$  $=$
$\big( {\sf wp}_{A \cup \{c\}\!}(G) \big)^{\rm rev}$, then we also have 
$(abc)^{\rm rev} = cba =_G \e$, i.e., $b^{-1} a^{-1} ba =_G \e$,
which implies $ba = ab$.
Hence, all generators in $A$ commute; it follows that $G$ is commutative.
 \ \ \  \ \ \ $\Box$

\bigskip

\noindent {\bf Examples.}

\smallskip

\noindent {\small (1)} The one-relator group with group presentation
$\,G = \langle \{a,b,c\} : \{abc\} \rangle$ has a word problem over 
$A = \{a,a^{-1},b,b^{-1},c,c^{-1}\}$ that is {\sl not closed under 
reversal}. If ${\sf wp}_{A\!}(G)$ were closed under reversal, then $G$
would be commutative, as we saw in the proof of Prop.\ \ref{PROPnotrev}.
The generators $c, c^{-1}$ and the relation can be eliminated by Tietze 
transformations, so $G$ is isomorphic to ${\sf FG}_2$, which is not
commutative.

\smallskip

\noindent {\small (2)} The free group ${\sf FG}_2$ with monoid generating
set $\{a, b, a^{-1}, b^{-1}\}$ has a word problem that is closed under 
reversal, as can be seen from its context-free grammar. But it is not 
generated by involutions, as the only involution in ${\sf FG}_2$ is the 
identity. (I thank the referee for this example.)

\smallskip

\noindent {\small (3)}  The {\em dihedral group} $D_{2n} = $
$\langle \{a,b\} : \{a^n, b^2, abab^{-1}\} \rangle$, with $|D_{2n}| = 2n$,
is non-commutative when $n \ge 3$.  
One can check that the word problem of $D_{2n}$ over $\{a,b\}$ is closed 
under reversal. The proof of Prop.\ \ref{PROPnotrev} yields a different
generating set for $D_{2n}$ for which the word problem is not closed under
reversal.  

\smallskip

\noindent {\small (4)} The {\em Thompson group} $V$ with the Higman 
generating set $\Gamma_{\! H}$ has a word problem that is closed under
reversal; see Cor.\ \ref{CORwpVrev} below. More 
generally, any non-commutative group that has a finite {\em generating set 
consisting of involutions} is an example; see Prop.\ \ref{PROPwpINVOLUTrev}. 
This includes many finite (simple) groups, and Coxeter groups.

\begin{pro} \label{PROPwpINVOLUTrev}
 \ Suppose $G$ is a group that has a finite generating set $A$ consisting of
{\em involutions}. Since we use internal generators, $\alpha$ and 
$\alpha^{-1}$ are the same letter in $A$, for all $\alpha \in A$. 
 \ Then we have 

\smallskip

\hspace{1,0cm} 
${\sf wp}_{A\!}(G) \,=\, \big({\sf wp}_{A\!}(G)\big)^{\rm rev}$.
\end{pro}
{\sc Proof.} For any $w = a_n \ldots a_1 \in A^*$ we have: $w =_G \e\,$ 
iff $\,\e =_G w^{-1} = a_1^{-1} \ldots a_n^{-1} = a_1 \ldots a_n$  $=$
$w^{\rm rev}$; the last equalities follows from $a_i = a_i^{-1}$.  
 \ \ \  \ \ \ $\Box$

\begin{cor} \label{CORwpVrev}
 \ The Thompson group $V$ has a finite generating set $\,\Gamma_{\!H}$
consisting of {\em involutions}, hence

\medskip

\hspace{1,0cm} ${\sf wp}_{\Gamma_{\! H}\!}(V)$  $\,=\,$ 
$\big({\sf wp}_{\Gamma_{\! H}\!}(V)\big)^{\rm rev}$.
\end{cor}
{\sc Proof.} We use the set $\Gamma_{\! H}$ of {\em Higman generators} for 
$V\,$ \cite[{\small p.\ 49}]{Hig74}, given by the tables

\medskip

{\sc not} $=$
$\left[ \hspace{-.08in} \begin{array}{l|l}
0 & 1  \\
1 & 0
\end{array} \hspace{-.08in} \right]$,
 \ \ \ $(1|01) =$ 
$\left[ \hspace{-.08in} \begin{array}{r|r|r}
00 & 01 & 1  \\
00 &  1 & 01
\end{array} \hspace{-.08in} \right]$,
 \ \ \ $(0|10) =$
$\left[ \hspace{-.08in} \begin{array}{r|r|r}
0 & 10 & 11  \\
10 & 0 & 11
\end{array} \hspace{-.08in} \right]$,
 \ \ \ $\tau_{1,2}$ $=$
$\left[ \hspace{-.08in} \begin{array}{r|r|r|r}
00 & 01 & 10 & 11  \\
00 & 10 & 01 & 11
\end{array} \hspace{-.08in} \right]$.

\medskip

\noindent Clearly, $a = a^{-1}$ for every $a \in \Gamma_{\! H}$. 
For more background information on $V$ and these generators, see the 
Subsection 2.3.  \ \ \   \ \ \  $\Box$

\bigskip

\noindent Not every finitely generated group is generated by a finite set
of involutions, hence one can ask:

\smallskip

\noindent {\bf Question:} Does every finitely generated group $G$ have 
some finite generating set $A$ such that ${\sf wp}_{A\!}(G)$ is closed
under reversal?

\subsection{Some explanations about the Thompson group $V$ and the 
Higman generators}

We give a brief introduction to the Thompson group $V$. This group can be 
defined in several rather different ways. Richard Thompson \cite{Th0, Th}
starts out with a few special permutations of the Cantor space
$\{0,1\}^{\omega}$ and lets $V$ be the group generated; he proves that $V$
is 2-generated, finitely presented and simple. 
In \cite{CFP} $V$ is a certain group of piecewise-linear bijections of the 
interval $[0,1] \cap \mathbb{Q}_2\,$ (where $\mathbb{Q}_2$ denotes the set
of binary rational numbers, i.e., numbers that can be expressed as $a/2^b$, 
with $a \in \mathbb{Z},\, b \in \mathbb{N}$).  
In \cite{Hig74} Higman describes $V$ as a certain group of isomorphisms 
of a free J\'onsson-Tarski algebra. This is formally the same as an action 
by permutations on the Cantor space and on $\{0,1\}^*$, and this is how $V$ 
is described by Scott \cite{EScott}. 
The same action as in \cite{EScott} is given in \cite{BiThomps}, but it is
now based on the terminology of finite maximal prefix codes, right ideals 
of $\{0,1\}^*$, and right-ideal morphisms or prefix replacements. 
The group $V$ also acts faithfully on $\{0,1\}^*\,0^{\omega}$. 

In this paper we use the description of \cite{BiThomps}, which formally 
looks the same as the ones in \cite{Hig74} and \cite{EScott}.
Every element of $V$ is a permutation $f$ of $\{0,1\}^{\omega}$ given by a 
{\em finite prefix-replacement table}  

\smallskip

\hspace{1,0cm}
$T \,=\, \left[ \hspace{-.08in} \begin{array}{r|r|r}
x_1 & \ldots & x_n  \\
y_1 & \ldots & y_n
\end{array} \hspace{-.08in} \right]$,

\smallskip

\noindent for any $n \ge 1$; equivalently, 
$\,T = \{(x_i,y_i) : 1 \le i \le n\}$ $\,\subseteq\,$ 
$\{0,1\}^* \x \{0,1\}^*$. Here 
$\,{\rm domC}(T) = \{x_i : 1 \le i \le n\}\,$ is called the domain code, 
and $\,{\rm imC}(T) = \{y_i : 1 \le i \le n\}\,$ is called the image code.  
Both ${\rm domC}(T)$ and ${\rm imC}(T)$ are chosen to be arbitrary finite 
maximal prefix codes of $\{0,1\}^*$ with the same cardinality. The 
correspondence $\,x_i \mapsto y_i\,$ (for $1 \le i \le n$) in the table is
a bijection from ${\rm domC}(T)$ onto ${\rm imC}(T)$. This is a 
recoding map, i.e., a translation from one code to another code.
The permutation $f$ of $\{0,1\}^{\omega}$ (or of $\{0,1\}^*\,0^{\omega}$)
is then defined by 

\smallskip

\hspace{1,0cm} $f(x_i\,z) \,=\, y_i\,z$ 

\smallskip

\noindent for every infinite bitstring $z \in \{0,1\}^{\omega}$.
Since ${\rm domC}(T)$ and ${\rm imC}(T)$ are finite maximal prefix codes, 
every element of $\{0,1\}^{\omega}$ can be written in a unique way as 
$x_i u$ and in a unique way as $y_i v$, for some $x_i \in$
${\rm domC}(T)$, $\,u \in \{0,1\}^{\omega}$, respectively for some 
$y_i \in {\rm imC}(T)$, $\,v \in \{0,1\}^{\omega}$. 
For any prefix code $P \subseteq \{0,1\}^*$ we have: $P\,\{0,1\}^{\omega}$
$=$ $\{0,1\}^{\omega}\,$ iff $\,P$ is a finite maximal prefix code. So, each
prefix-replacement table defines a permutation of $\{0,1\}^{\omega}$.

The prefix-replacement table also defines an injective partial function $f$
on $\{0,1\}^*$, by using the same formula $\,f(x_i z) = y_i z$, but with 
$z \in \{0,1\}^*$. Now $f(x_i) = y_i$, hence $\,f(x_i z) = f(x_i)\,z$; it 
follows that on $\{0,1\}^*$, $f$ is a right-ideal isomorphism from the 
right ideal $\,{\rm Dom}(T) = {\rm domC}(T)\,\{0,1\}^*$ onto the right ideal
$\,{\rm Im}(T) = {\rm imC}(T)\,\{0,1\}^*$. So the permutation $f$ of 
$\{0,1\}^{\omega}$ is extended to an injective partial function on 
$\,\{0,1\}^{\omega} \cup \{0,1\}^*$, also called $f$; 
the original $f$ is now the restriction $f|_{\{0,1\}^{\omega}}$.

It is a simple exercise to show that the composition of two permutations
given by finite prefix-replacement tables is also described by a finite 
prefix-replacement table.

The same permutation $f$ can be described by infinitely many different 
prefix-replacement tables. Indeed, for the above table 
$\,T = \{(x_i,y_i) : 1 \le i \le n\}\,$ for $f$, the table 

\smallskip

\hspace{1,0cm} $\,T_j \ = \ \big(T \minus \{(x_j,y_j)\}\big)$  $\,\cup\,$ 
$\{(x_j0, y_j0),\, (x_j1, y_j1)\}\,$ 

\smallskip

\noindent also defines $f\!: \{0,1\}^{\omega} \to \{0,1\}^{\omega}$; this 
holds for any $j$ with $1 \le j \le n$ (where $n$ is the cardinality of
$T$).  The passage from $T$ to $T_j$ is called the one-step 
{\em restriction} at $j$; indeed, the right ideals $\,{\rm Dom}(T_j)$ $=$
$\big(\{x_i : 1 \le i \le n,\,i \ne j\} \cup \{x_j0,x_j1\}\big)\,\{0,1\}^*$ 
and $\,{\rm Im}(T_j)$ $=$
$\big(\{y_i : 1 \le i \le n,\,i \ne j\} \cup \{y_j0,y_j1\}\big)\,\{0,1\}^*$, 
are strict subideals of the corresponding right ideals ${\rm Dom}(T)$ and 
${\rm Im}(T)$.  Likewise, the passage from $T_j$ to $T$ is called the 
one-step {\em extension} of $T_j$ at $j$; extension steps make the ideals
larger. And there exists a unique maximal extension for $f$, as a 
right-ideal morphism of $\{0,1\}^*$.

Any table $T$ with $x_i = y_i$ for all $i$, represents the identity 
function on $\{0,1\}^{\omega}$; conversely, all tables of the identity are
of that form. By applying extensions, this table can be transformed to 
$\{(\e,\e)\}$. For the inverse table $\,T^{-1}$  $=$ 
$\{(y_i,x_i) : 1 \le i \le n\}$, the composites $T \circ T^{-1}$ and 
$T^{-1} \circ T$ represent the identity, so $V$ is a group. The 
representation of $V$ by prefix-replacement tables makes it easy to solve 
the word problem of $V$; the obvious procedure is in ${\sf DTime}(n^3)$.

\smallskip

Higman did not give meaningful names to his generators. Our names are
motivated as follow: {\sc not} gives the logical negation in the first
bit of an input string in $\{0,1\}^{\omega}$; 
$(1|01)$ swaps the prefixes 1 and 01 of an input string; 
$(0,10)$ swaps 0 and 10; and $\tau_{1,2}$ transposes the positions of the 
first two bits of the input string. 
Higman \cite{Hig74} gives an explicit finite presentation of $V$ over
these generators. Besides being involutions, they have the property
that their tables only contain bitstrings of length $\,\le 2$. 

\medskip

\section{The deterministic complexity of the word problem of $V$}

\subsection{Proof of the main result}

For any prefix-replacement table $T$, and for any finite set $S$ of
tables, we define

\smallskip

\hspace{1,0cm} ${\rm maxlen}(T) \ = \  \max \{\,|z| \,:\,$
$z \in {\rm domC}(T) \,\cup\, {\rm imC}(T) \,\}$,

\smallskip

\hspace{1,0cm} ${\rm maxlen}(S) \ = \ $ 
$\max \{\,{\rm maxlen}(T) \,:\, T \in S\,\}$.

\medskip

\noindent
Let $A$ be a finite monoid generating set of the Thompson group $V$. 
The table obtained by composing the generators in $w \in A^*$ is denoted by 
$w(.)$. See the definition of $V$ in Subsection 2.3. By 
\cite[{\small Cor.\ 3.7}]{BiThomps}: 

\smallskip

\hspace{1,0cm} ${\rm maxlen}(w(.)) \ \le \ |w| \ {\rm maxlen}(A)$.  

\smallskip

\noindent For $x \in \{0,1\}^*$ let $w(x)$ be the result of applying the 
prefix replacement, given by the table $w(.)$, to $x$. Then we have: 

\smallskip

\hspace{1,0cm} If $x \in \{0,1\}^*$ satisfies 
$ \ |x| \,\ge\, |w| \ {\rm maxlen}(A) \ $ then 
$w(x)$ is defined, and 

\hspace{1,0cm} $w(x)$ can be computed by successively applying the 
generators in $w$ 

\hspace{1,0cm} (without applying maximum extensions).

\smallskip

\noindent For the action of the elements of $V$ as injective partial 
functions on $\{0,1\}^*$ we have:  for all $w \in A^*$, 

\smallskip

\hspace{1,0cm} $w \in {\sf cowp}_A(V)$ \ iff \ there exists 
$x \in \{0,1\}^*$ such that $w(x)$ is defined and $w(x) \ne x$.

\smallskip

\noindent This holds because the identity of $V$ is represented by any 
prefix-replacement table with two equal rows (see Subsection 2.3).

\smallskip

The following Lemma plays a crucial role in Lehnert and Schweitzer's proof
that the word problem of $V$ is in {\sf coCFL}; the Lemma is intuitive and 
does not appear explicitly in \cite{LehSch}. 

\smallskip

\begin{lem} \label{LEMcompdepthV} {\bf (narrow point).} 
 \ Let $\,B = B^{-1} \subseteq V$ be a finite subset, and let 
$S = \langle B \rangle$ be the subgroup of $V$ generated by $B$. 

\smallskip

\noindent For any $\,w = b_n \ldots b_1 \in B^+\,$ (with $b_i \in B$ for
$1 \le i \le n$), and any $\,x \in  \{0,1\}^*$, let

\medskip

\hspace{1,0cm} $x = x_0$  $\stackrel{b_1}{\longmapsto}$ 
$x_1$  $\stackrel{b_2}{\longmapsto}$ \ $\ldots$ 
 \ $\stackrel{b_{i-1}}{\longmapsto}$
$x_{i-1}$   $\stackrel{b_i}{\longmapsto}$ 
$x_i$   $\stackrel{b_{i+1}}{\longmapsto}$ \ $\ldots$
 \ $\stackrel{b_n}{\longmapsto}$
$x_n = w(x)$

\medskip

\noindent be the action steps of $w$ on input $x = x_0$, as a 
composition of prefix replacement tables, where 
$\,x_i = b_i(x_{i-1})\,$ for $1 \le i \le n$.  
We assume that $w(x)$ is defined, i.e., $b_i(x_{i-1})$ is defined for
all $i$.

Note the order $\,b_n \ldots b_1$, since $V$ and $S$ act on the {\em left} by
$\,w(x) \,=\, b_n( \ \ldots \ b_2(b_1(x))...)\,$.

\medskip

\noindent Then there exist $\,s, z_0, z_1, \,\ldots\,, z_n \in \{0,1\}^*$
such that

\medskip

\noindent {\small \bf (1)} \ \ $x_i = z_i \,s$ \ \ and 
 \ \ $|z_i| \le |w| \ {\rm maxlen}(B)$, \ \ for all 
$i = 0,1,\,\ldots,\,n$.

\medskip

\noindent {\small \bf (2)} \ \ The following are the action steps of $w(.)$
on input $z_0:$ 

\medskip

\hspace{1,5cm} $z_0$  $\stackrel{b_1}{\longmapsto}$
$z_1$  $\stackrel{b_2}{\longmapsto}$ \ \ $\ldots$
 \ \ $\stackrel{b_{i-1}}{\longmapsto}$
$z_{i-1}$   $\stackrel{b_i}{\longmapsto}$
$z_i$   $\stackrel{b_{i+1}}{\longmapsto}$ \ \ $\ldots$
 \ \ $\stackrel{b_n}{\longmapsto}$
$z_n = w(z_0)$,

\medskip

 \ \ \ such that $b_i(z_{i-1})$ is defined and $\,z_i = b_i(z_{i-1})\,$ 
for all $i = 1,\,\ldots,\,n$.

\medskip

\noindent {\small \bf (3)} \ \ There exists $\,k \in$ 
$\{0,1, \ldots, n\}$ such that \ $|z_k| \le {\rm maxlen}(B)$.  The
position $k$ in $w$ is called a {\em narrow point}, since $z_k$ has 
bounded length (with a bound that does not depend on $w$).
See {\rm Fig.\ 1}.
 
\medskip

\noindent {\small \bf (4)} \ \ $w \in {\sf cowp}_{B\!}(S)$ \ iff
 \ there exists $z_0 \in \{0,1\}^*$ such that $w(z_0) \ne z_0$, and the 
sequence $\,z_i = b_i \ldots b_1(z_0)\,$ (for $i \in [0,n]$) satisfies 
{\small \rm (2)} and {\small \rm (3)}. \ Moreover,

\smallskip

\hspace{0,5cm} $w \in {\sf cowp}_{B\!}(S)$ \ iff 

\hspace{0,5cm} there exists $z \in \{0,1\}^*$ such that 
$\,|z| \,\le\, {\rm maxlen}(B)$, and 
$\,b_k \ldots b_1 b_n \ldots b_{k+1}(z 0^{\omega})$ $\ne$ $z 0^{\omega}$. 

\hspace{0,5cm} (If $k = 0$ then $\,b_k \ldots b_1 = \e$; if $k = n$
then $b_n \ldots b_{k+1} = \e$.)
\end{lem}
{\sc Proof.} (1), (2)  When an element $b_i \in V$ is applied to a word 
$x_i \in {\rm Dom}(b_i)$, a prefix of $x_i$ of length $\,\le$
${\rm maxlen}(b_i)\,$ is modified (see the explanations before Lemma
\ref{LEMcompdepthV}). So after $n = |w|\,$ steps, a prefix of $x_i$ of
length $\,\le n \, {\rm maxlen}(B)\,$ has been modified. This implies 
items (1) and (2).

\smallskip

\noindent (3) By choosing the prefixes $z_i$ of $x_i$ to have minimal 
length, subject to (1) and (2) we obtain $|z_k| \le {\rm maxlen}(B)$, for 
some $k$. Indeed, if {\em all} $z_i$ were such that 
$|z_i| > {\rm maxlen}(B)$, then every step
$\,z_i \longmapsto b_{i+1}(z_i) = z_{i+1}\,$ 
would change only a strict prefix of $z_i$ (for all $i = 0,1,\ldots,n-1$). 
So $s$ could be lengthened, and every $z_i$ could be shortened, while (1) 
and (2) would still hold.  See Figure 1. 

There could be more than one location $k$ in
$\,b_n\,\dots\,b_1\,$ at which $\,|z_k|\,$ reaches the minimum; those
minimum locations depend on the input $z_0$.

\bigskip

\unitlength=.8mm
\begin{picture}(80,95)

\put(15,10){\line(0,1){70}}         
\put(15,40){\makebox(0,0)[cc]{-}}
\put(15,80){\makebox(0,0)[cc]{-}}
\put(12,25){\makebox(0,0)[cc]{$s$}}
\put(12,60){\makebox(0,0)[cc]{$z_0$}}

\put(25,10){\line(0,1){75}}         
\put(25,40){\makebox(0,0)[cc]{-}}
\put(25,85){\makebox(0,0)[cc]{-}}

\put(35,10){\line(0,1){80}}
\put(35,40){\makebox(0,0)[cc]{-}}
\put(35,90){\makebox(0,0)[cc]{-}}

\put(45,10){\line(0,1){75}}
\put(45,40){\makebox(0,0)[cc]{-}}
\put(45,85){\makebox(0,0)[cc]{-}}

\put(55,10){\line(0,1){65}}
\put(55,40){\makebox(0,0)[cc]{-}}
\put(55,75){\makebox(0,0)[cc]{-}}

\put(65,10){\line(0,1){55}}
\put(65,40){\makebox(0,0)[cc]{-}}
\put(65,65){\makebox(0,0)[cc]{-}}

\put(75,10){\line(0,1){50}}
\put(75,40){\makebox(0,0)[cc]{-}}
\put(75,60){\makebox(0,0)[cc]{-}}

\put(85,10){\line(0,1){40}}
\put(85,40){\makebox(0,0)[cc]{-}}
\put(85,50){\makebox(0,0)[cc]{-}}

\put(95,10){\line(0,1){35}}         
\put(95,40){\makebox(0,0)[cc]{-}}
\put(95,45){\makebox(0,0)[cc]{-}}
\put(92,25){\makebox(0,0)[cc]{$s$}}
\put(92,43){\makebox(0,0)[cc]{$z_k$}}

\put(105,10){\line(0,1){40}}
\put(105,40){\makebox(0,0)[cc]{-}}
\put(105,50){\makebox(0,0)[cc]{-}}

\put(115,10){\line(0,1){45}}
\put(115,40){\makebox(0,0)[cc]{-}}
\put(115,55){\makebox(0,0)[cc]{-}}

\put(125,10){\line(0,1){50}}
\put(125,40){\makebox(0,0)[cc]{-}}
\put(125,60){\makebox(0,0)[cc]{-}}

\put(135,10){\line(0,1){60}}
\put(135,40){\makebox(0,0)[cc]{-}}
\put(135,70){\makebox(0,0)[cc]{-}}

\put(145,10){\line(0,1){65}}
\put(145,40){\makebox(0,0)[cc]{-}}
\put(145,75){\makebox(0,0)[cc]{-}}

\put(155,10){\line(0,1){70}}
\put(155,40){\makebox(0,0)[cc]{-}}
\put(155,80){\makebox(0,0)[cc]{-}}

\put(165,10){\line(0,1){65}}
\put(165,40){\makebox(0,0)[cc]{-}}
\put(165,75){\makebox(0,0)[cc]{-}}

\put(175,10){\line(0,1){65}}            
\put(175,40){\makebox(0,0)[cc]{-}}
\put(175,75){\makebox(0,0)[cc]{-}}
\put(172,25){\makebox(0,0)[cc]{$s$}}
\put(172,57){\makebox(0,0)[cc]{$z_n$}}

\end{picture}

\vspace{-3mm}

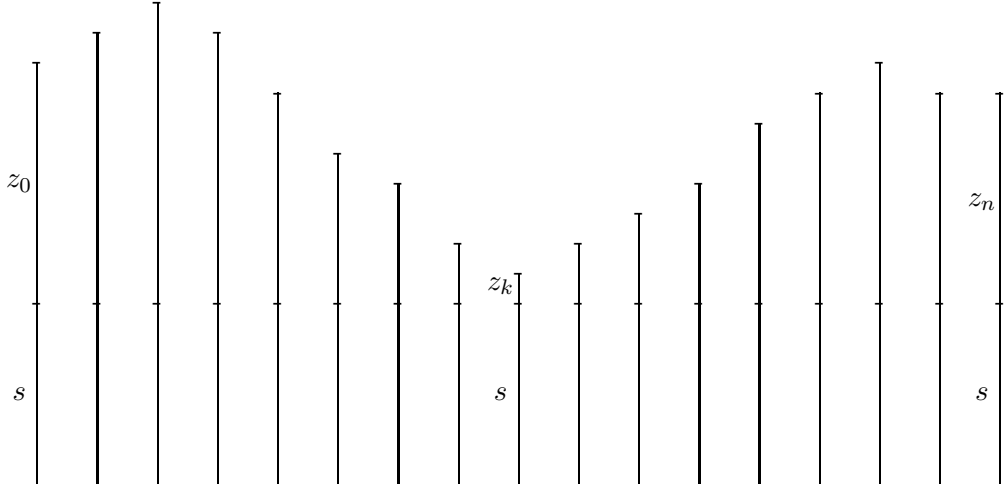
\noindent {\small {\bf Figure 1:} \ Illustration (with $n = 16$) of a 
computation \ $b_n \ldots b_1()$: 
$\,z_0 s \in \{0,1\}^* \longmapsto z_n s \in \{0,1\}^*\,$ in $V,$ with 
minimum length at $z_k s$, with $|z_k| \le {\rm maxlen}(B)$.} 

\bigskip

\medskip

\noindent {\sc Proof} of Lemma \ref{LEMcompdepthV} cont'd: 
(4) $\,w \in {\sf cowp}_{B\!}(S)\,$ iff there exists
$z_0 \in \{0,1\}^*$ such that $w(z_0) \ne z_0$, and such that the strings
$z_i = b_i \ldots b_1(z_0)$ are obtained as in (1) and satisfy (2), and (3).
I.e., for all $i \in [0,n]$,
$\,b_n \ldots b_{i+1}(z_i)$ $\ne$ $z_0$. 

Hence we also have $\,b_i \ldots b_1 b_n \ldots b_{i+1}(z_i 0^{\omega})$
$\ne$ $z_i 0^{\omega} = b_i \ldots b_1(z_0 0^{\omega})$ for all 
$i \in [0,n]$, since every generator is invertible and is defined on all of
$\{0,1\}^*\,0^{\omega}$.
By {\small (3)} there exists $k \in [0,n]$ such that 
$\,|z_k| \le {\rm maxlen}(B)$.
And $\,b_k \ldots b_1 b_n \ldots b_{k+1}(z_k 0^{\omega})$ $\ne$
$z_k 0^{\omega}$, since this holds for all $i$.
In turn, $\,b_k \ldots b_1 b_n \ldots b_{k+1}(z_k 0^{\omega})$ $\ne$
$z_k 0^{\omega}\,$ implies $\,w \in {\sf cowp}_{B\!}(S)$, by closure of
${\sf cowp}_{B\!}(S)$ under cyclic permutation.
   \ \ \  \ \ \ $\Box$

\medskip

\noindent {\bf Remark.} In (4) we use the faithful action of $V$ on
$\{0,1\}^*\,0^{\omega}$; then 
$b_i \ldots b_1 b_n \ldots b_{i+1}(z_i 0^{\omega})$ is defined
for all $i \in [0,n]$. On the other hand, although $w(z_0)$ is defined,
$b_i \ldots b_1 b_n \ldots b_{i+1}(z_i)$ $=$ $b_i \ldots b_1 w(z_0)$
might be undefined for the partial action on $\{0,1\}^*$.

\bigskip

\medskip

\noindent {\bf Remark about the Brin-Thompson group $2V$:}
Item (3) of Lemma \ref{LEMcompdepthV} is actually subtle, as an attempt to 
apply it to $2V$ shows. We follow \cite{JCBnG} for the definition of $2V$.
Let $A_2$ be a finite monoid generating set of the group $2V$, and let us 
define length by $\,|x|$  $=$ $\max\{|x^{(1)}|,\, |x^{(2)}|\}\,$ for any 
$x = (x^{(1)},x^{(2)})$  $\in$ $\{0,1\}^* \x \{0,1\}^*$. 
Then the reasoning in the proof of Lemma \ref{LEMcompdepthV} seems to work, 
at first look. But this would lead (via a proof similar to the one for 
Prop.\ \ref{PROwpVcyc4DCF}) to the conclusion that the word problem of $2V$ 
is in {\sf P} (while it is also {\sf coNP}-complete by \cite{JCBnG}).  
In fact Lemma \ref{LEMcompdepthV}(3) does not hold in $2V$, because the 
minimum of $|z_i^{(1)}|\,$ (1st coordinate) can be in a different location
in $a_n \,\ldots\, a_1$ than the minimum of $|z_j^{(2)}|\,$ (2nd 
coordinate); this is illustrated in Figure 2. This can also be seen in the 
example of the shift $\sigma \in 2V$, where $\,\sigma^n(.)$: 
$(\e, u) \mapsto (u^{\rm rev}, \e)\,$ for any word $u \in \{0,1\}^+$ with 
$|u| = n$.

\unitlength=.8mm
\begin{picture}(80,100)

\put(15,10){\line(0,1){70}}         
\put(15,40){\makebox(0,0)[cc]{-}}
\put(15,80){\makebox(0,0)[cc]{-}}
\put(16,6){\makebox(0,0)[cc]{$s$}}
\put(16,82){\makebox(0,0)[cc]{$z_0$}}
\put(18,10){\line(0,1){50}}         
\put(18,40){\makebox(0,0)[cc]{-}}
\put(18,60){\makebox(0,0)[cc]{-}}

\put(25,10){\line(0,1){75}}         
\put(25,40){\makebox(0,0)[cc]{-}}
\put(25,85){\makebox(0,0)[cc]{-}}
\put(28,10){\line(0,1){55}}         
\put(28,40){\makebox(0,0)[cc]{-}}
\put(28,65){\makebox(0,0)[cc]{-}}

\put(35,10){\line(0,1){80}}
\put(35,40){\makebox(0,0)[cc]{-}}
\put(35,90){\makebox(0,0)[cc]{-}}
\put(38,10){\line(0,1){55}}
\put(38,40){\makebox(0,0)[cc]{-}}
\put(38,65){\makebox(0,0)[cc]{-}}

\put(45,10){\line(0,1){75}}
\put(45,40){\makebox(0,0)[cc]{-}}
\put(45,85){\makebox(0,0)[cc]{-}}
\put(48,10){\line(0,1){60}}
\put(48,40){\makebox(0,0)[cc]{-}}
\put(48,70){\makebox(0,0)[cc]{-}}

\put(55,10){\line(0,1){65}}
\put(55,40){\makebox(0,0)[cc]{-}}
\put(55,75){\makebox(0,0)[cc]{-}}
\put(58,10){\line(0,1){68}}
\put(58,40){\makebox(0,0)[cc]{-}}
\put(58,78){\makebox(0,0)[cc]{-}}

\put(65,10){\line(0,1){55}}
\put(65,40){\makebox(0,0)[cc]{-}}
\put(65,65){\makebox(0,0)[cc]{-}}
\put(68,10){\line(0,1){70}}
\put(68,40){\makebox(0,0)[cc]{-}}
\put(68,80){\makebox(0,0)[cc]{-}}

\put(75,10){\line(0,1){50}}
\put(75,40){\makebox(0,0)[cc]{-}}
\put(75,60){\makebox(0,0)[cc]{-}}
\put(78,10){\line(0,1){75}}
\put(78,40){\makebox(0,0)[cc]{-}}
\put(78,85){\makebox(0,0)[cc]{-}}

\put(85,10){\line(0,1){40}}
\put(85,40){\makebox(0,0)[cc]{-}}
\put(85,50){\makebox(0,0)[cc]{-}}
\put(88,10){\line(0,1){70}}
\put(88,40){\makebox(0,0)[cc]{-}}
\put(88,80){\makebox(0,0)[cc]{-}}

\put(95,10){\line(0,1){35}}         
\put(95,40){\makebox(0,0)[cc]{-}}
\put(94.8,40){\line(0,1){5}}
\put(95.2,40){\line(0,1){5}}
\put(95,45){\makebox(0,0)[cc]{-}}
\put(96,6){\makebox(0,0)[cc]{$s$}}
\put(97,78){\makebox(0,0)[cc]{$z_i$}}
\put(98,10){\line(0,1){65}}         
\put(98,40){\makebox(0,0)[cc]{-}}
\put(98,75){\makebox(0,0)[cc]{-}}

\put(105,10){\line(0,1){40}}
\put(105,40){\makebox(0,0)[cc]{-}}
\put(105,50){\makebox(0,0)[cc]{-}}
\put(108,10){\line(0,1){60}}
\put(108,40){\makebox(0,0)[cc]{-}}
\put(108,70){\makebox(0,0)[cc]{-}}

\put(115,10){\line(0,1){45}}
\put(115,40){\makebox(0,0)[cc]{-}}
\put(115,55){\makebox(0,0)[cc]{-}}
\put(118,10){\line(0,1){55}}
\put(118,40){\makebox(0,0)[cc]{-}}
\put(118,65){\makebox(0,0)[cc]{-}}

\put(125,10){\line(0,1){50}}
\put(125,40){\makebox(0,0)[cc]{-}}
\put(125,60){\makebox(0,0)[cc]{-}}
\put(128,10){\line(0,1){50}}
\put(128,40){\makebox(0,0)[cc]{-}}
\put(128,60){\makebox(0,0)[cc]{-}}

\put(135,10){\line(0,1){60}}
\put(135,40){\makebox(0,0)[cc]{-}}
\put(135,70){\makebox(0,0)[cc]{-}}
\put(138,10){\line(0,1){40}}
\put(138,40){\makebox(0,0)[cc]{-}}
\put(138,50){\makebox(0,0)[cc]{-}}

\put(145,10){\line(0,1){65}}           
\put(145,40){\makebox(0,0)[cc]{-}}
\put(145,75){\makebox(0,0)[cc]{-}}
\put(148,10){\line(0,1){35}}
\put(148,40){\makebox(0,0)[cc]{-}}
\put(148,45){\makebox(0,0)[cc]{-}}
\put(147.8,40){\line(0,1){5}}
\put(148.2,40){\line(0,1){5}}
\put(146,6){\makebox(0,0)[cc]{$s$}}
\put(147,77){\makebox(0,0)[cc]{$z_j$}}

\put(155,10){\line(0,1){70}}
\put(155,40){\makebox(0,0)[cc]{-}}
\put(155,80){\makebox(0,0)[cc]{-}}
\put(158,10){\line(0,1){45}}
\put(158,40){\makebox(0,0)[cc]{-}}
\put(158,55){\makebox(0,0)[cc]{-}}

\put(165,10){\line(0,1){65}}
\put(165,40){\makebox(0,0)[cc]{-}}
\put(165,75){\makebox(0,0)[cc]{-}}
\put(168,10){\line(0,1){55}}
\put(168,40){\makebox(0,0)[cc]{-}}
\put(168,65){\makebox(0,0)[cc]{-}}

\put(175,10){\line(0,1){65}}            
\put(175,40){\makebox(0,0)[cc]{-}}
\put(175,75){\makebox(0,0)[cc]{-}}
\put(177,78){\makebox(0,0)[cc]{$z_n$}}
\put(177,6){\makebox(0,0)[cc]{$s$}}
\put(178,10){\line(0,1){60}}            
\put(178,40){\makebox(0,0)[cc]{-}}
\put(178,70){\makebox(0,0)[cc]{-}}

\end{picture}

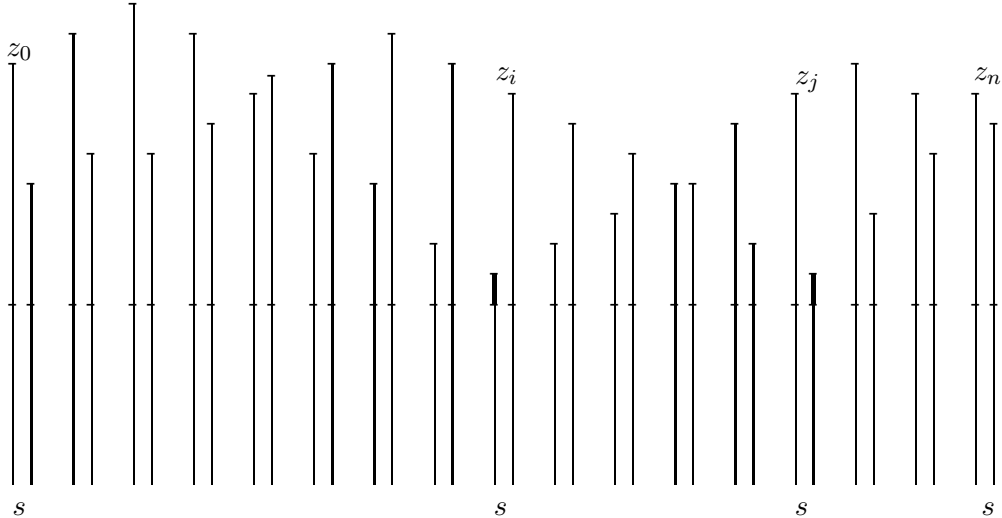
\noindent {\small {\bf Figure 2:} Illustration (with $n = 16$) of a 
computation

\hspace{1,0cm} $a_n \ldots a_1()$: $\,z_0 s \in \{0,1\}^* \x \{0,1\}^*$
$\longmapsto$  $z_n s \in \{0,1\}^* \x \{0,1\}^*$ 

\noindent  in $2V$, with minimum length in coordinate 1 at $z_i s$, and 
minimum length in coordinate 2 at $z_j s$, where 
$\,|z_i^{(1)}|,\, |z_j^{(2)}| \le {\rm maxlen}(A_2)$.}

\bigskip

\medskip

\begin{lem} \label{LEMz0omega}
 \ For any $s, z \in \{0,1\}^+$ such that 
$\,|z| \le |s|$ we have:

\medskip

$s 0^{\omega} \ne z 0^{\omega}$ \ \ \ iff
  \ \ \ $\,s \ne z 0^{|s| - |z|}$ 

\smallskip

iff \ \ $z$ is not a prefix of $s$, \ or
 \ $s = z t\,$ for some $\,t \in \{0,1\}^+ \minus 0^+$ $=$ 
$\{0,1\}^*\,1\, \{0,1\}^*$.
\end{lem}
{\sc Proof.} In general, $x_1, x_2 \in \{0,1\}^{\omega}$ are different iff 
there exists $m \ge 1$ such that the respective prefixes of $x_1$ and $x_2$ 
of length $m$ are different.
If $x_1 = z 0^{\omega}$ and $x_2 = s 0^{\omega}$ for some $s, z \in$
$\{0,1\}^+$, and $|z| \le |s|$, then $x_1 \ne x_2$ iff there is a difference
in the non-$0^{\omega}$ part, i.e., iff $\,s \ne z 0^{|s| - |z|}$. 
If $|z| \le |s|$ then $s = z' t$ for some $z', t \in \{0,1\}^*$ with 
$|z'| = |z|$, and $|t| = |s| - |z|$. Then $\,s \ne z 0^{|s| - |z|}$ iff
$z \ne z'$ (i.e., $z$ is not a prefix of $s$), or $t \ne 0^{|s| - |z|}$ 
(i.e., $t$ contains the letter $1$). 
 \ \ \  \ \ \ $\Box$

\bigskip

\noindent
In the next Lemma we construct a deterministic push-down automaton (dpda).
The Appendix describes the notation used here. 
The Lemma is a stronger form of \cite[{\small Sect.\ 5, Step 2}]{LehSch}, 
where a nondeterministic pushdown automaton (pda) was constructed, while
$(.)^{\rm rev}$ was omitted (since {\sf CFL} is closed under reversal). 
Moreover, the Lemma holds for any finitely generated subgroup of $V$.

\begin{lem} \label{LEMdcflL_Z} 
 \ Let $B = B^{-1}$ be a finite subset of $V,$ generating a subgroup 
$\,\langle B \rangle \,\subseteq\, V$.  For any $z \in \{0,1\}^+$ let

\medskip

\hspace{1,5cm} $L_z \,=\, \{\,w \in B^+ :\,$
   $w(z 0^{\omega}) \ne z 0^{\omega}\,\}$
 \ \ \ $(\,\subseteq\,$ ${\sf cowp}_B(\langle B \rangle))$.

\medskip

\noindent Then $L_z^{\rm \,rev}$ is {\em deterministic} context-free.
\end{lem}
{\sc Proof.} For any $w \in L_z$ we have $w(.) \ne {\sf id}(.)$, so 
$L_z \subseteq {\sf cowp}_B(V)$. 

The reason why $L_z^{\rm \,rev}$ is in {\sf DCFL} (and not $L_z$ itself) is 
that in $w(.) = a_n \ldots a_1(.)$, the functions $a_i$ are applied from 
right to left, since $V$ acts on $\{0,1\}^*$ and $\{0,1\}^{\omega}$ on the 
{\sl left}.  Below, the action of $V$ will be carried out by a deterministic
push-down automaton (dpda), which reads the letters in the order in which 
$V$ acts on $\{0,1\}^*$, i.e., starting with $a_1$ and ending with $a_n$.  

\smallskip

\noindent {\sf Definition and Remarks:}
An {\em endmarker language} consists of an alphabet $B$, a letter 
$\# \not\in B$ (called {\em input endmarker}), and a subset of 
$B^* \#\,$; so the language has the form $\,L \#\,$ over the 
alphabet $B \cup \{\#\}$, with $L \subseteq B^*$.  It was proved in 
\cite{GinsGrei66} that $L \#$ is in {\sf DCFL} iff $L$ is in {\sf DCFL};
see also e.g.\ \cite[{\small Thm.\ 11.2.2 and 11.3.1}]{Harrison},
\cite[{\small Thm.\ 2.43}]{Sipser}, \cite[{\small Thm.\ 10.2}]{HU}.
\hspace{1,0cm} {\small [End, Remarks]}

\smallskip 
 
We now construct a deterministic push-down automaton (dpda) that accepts the
endmarker language $L_z^{\rm \,rev} \#$; then by the Definition and Remarks 
above, $L_z^{\rm \,rev}$ is in {\sf DCFL}.
See the Appendix for information about the (deterministic) pda{\it s}.
 
Our dpda has the state set $\{q_0, q_1, q_{\rm a}\}$. The input alphabet is 
$B \cup \{\#\}$, and the stack alphabet is $\{0,1,\bot\}$, where $\bot$ is 
the bottom marker of the stack. The start configuration is $(q_0, z \bot)$,
where $z \in \{0,1\}^+$ is the fixed string that defines $L_z$. 
The dpda accepts an input $w\# \,\in\, B^*\#\,$ iff the dpda reaches the 
accept state $q_{\rm a}$ after reading $w\#\,$ (acceptance by final state). 

A transition of the dpda applies the next input letter $b \in B$ to the part 
$x \in \{0,1\}^*$ of the current stack content $x\bot$.  The action in the 
Thompson group is a right-ideal morphism that does a prefix replacement, 
transforming $x$ into $b(x)$. This is exactly what a pda does on the stack: 
the Thompson group elements treat their argument like a stack, where a 
bitstring $x = x_k \ldots x_1$ is stored with $x_1$ towards the stack-bottom 
and $x_k$ at the top.  If $b(x)$ is undefined, we use the action of $V$ on
$\{0,1\}^*\,0^{\omega}$ and let the dpad treat $x$ as $x0^m$, where $m$ is 
exactly so that $x0^m \in {\rm domC}(b)$. Indeed, ${\rm domC}(b)$ is a 
finite maximal prefix code, hence for every $x$ that is too short to belong
to ${\rm Dom}(b)$, there is exactly one $m \in \mathbb{N}$ such that 
$x0^m \in {\rm domC}(b)$ (see \cite[{\small Section 1}]{BiThomps}). While
the dpada simulates the action of $V$, the stack content $x \bot$ is never 
just $\bot$ (with $x = \e$), since ${\rm imC}(b)$ does not contain $\e$. 

When the endmarker $\#$ is encountered in the input, the dpda starts a 
process that checks whether the current stack content 
$s \bot \in \{0,1\}^* \bot$ satisfies $s 0^{\omega} \ne z 0^{\omega}\,$ 
(where $z$ is the fixed bitstring that defines $L_z$). 
The dpda uses transitions that can read a certain depth into the stack (but 
by a bounded amount, since the set of transitions is finite and fixed).
On input $\#$ there are three cases, based on Lemma \ref{LEMz0omega}:

\smallskip

\noindent {\small (1)} $|s| \le |z|$: 
 \ Then $z 0^{\omega} \ne s 0^{\omega}$ iff $z \ne s 0^{|s| - |z|}$. The 
latter can be checked by looking into the stack $s\bot$ by at depth at
most $|z|$. If $z \ne s 0^{|s| - |z|}$ then the dpda goes to state 
$q_{\rm a}$ by a transition on input $\#$.
If $z = s 0^{|s| - |z|}$ then the dpda rejects by having no transition. 

\noindent {\small (2)} $|s| > |z|$:  
 \ Then $z 0^{\omega} \ne s 0^{\omega}$ iff $z$ is not a prefix of $s$, or
$s = zt$, for some $t \in \{0,1\}^* 1 \{0,1\}^*$.

\noindent {\small (2.1)}
If $z$ is not a prefix of $s$ (which can be checked by looking into the 
stack by a depth at most $|z|$), then the dpda goes to state $q_{\rm a}$
by a transition on input $\#$.

\noindent {\small (2.2)}
If $z$ is a prefix of $s = zt$, for some $t \in \{0,1\}^* 1 \{0,1\}^*$, the
dpda sees that $z$ is a strict prefix of $s$ by looking into the stack by
the bounded depth $|z| + 1$. To check whether $t$ contains at least one $1$,
the dpda pops $z$ on input $\#$ and goes to state $q_1$. In state $q_1$,
with $\e$-transitions, the dpda pops letters of $t$ one at a time until 
either a letter $1$ is found (it goes to state $q_{\rm a}$ in that
case), or $\bot$ is encountered (it rejects in that case; it has no 
transition $\,(q_1, \bot) \stackrel{\e}{\to} \ldots \,$).

\smallskip

In state $q_{\rm a}$ the dpda does not have any transitions, hence no word 
in $\,L_z \# B^+\,$ will be accepted. 

\medskip

\noindent In detail, the set of transitions of the dpda is

\medskip

 \ \ \ $\{\,(q_0, r) \stackrel{b}{\to} (q_0, b(r)) :\,$
$b \in B, \ r \in {\rm domC}(b) \,\}$           

\medskip

 \ \ \ $\,\cup\,$
 \ \ \ $\{\,(q_0, r \bot) \stackrel{b}{\to} (q_0, b(r 0^m)\,\bot) :\,$
$b \in B, \ r \in \{0,1\}^+, \ m \ge 1, \ r 0^m \in {\rm domC}(b) \,\}$

\medskip

 \ \ \ $\,\cup\,$ 
 \ \ \ $\{\,(q_0, s \bot) \stackrel{\#}{\to} (q_{\rm a}, s \bot) :\,$ 
$s \in \{0,1\}^+, \ |s| < |z|, \ z \ne s 0^{|z| - |s|} \,\}$        

\medskip

 \ \ \ $\,\cup\,$
 \ \ \ $\{\,(q_0, s) \stackrel{\#}{\to} (q_{\rm a}, s) :\,$
$s \in \{0,1\}^*, \ |s| = |z|, \ z \ne s \,\}$          

\medskip

 \ \ \ $\,\cup\,$
 \ \ \ $\{\,(q_0, za) \stackrel{\#}{\to} (q_1, a) :\, a \in \{0,1\}\,\}$

\medskip

 \ \ \ $\,\cup\,$
 \ \ \ $\{\,(q_1, 0) \stackrel{\e}{\to} (q_1, \e) \,\}$

\medskip

 \ \ \ $\,\cup\,$
 \ \ \ $\{\,(q_1, 1) \stackrel{\e}{\to} (q_{\rm a}, 1) \,\}$.

\medskip

\noindent The first two sets of transitions carry out a simulation of the 
action of $V$.  The remaining sets let the dpda recognize whether the stack 
content $s\bot$ satisfies $s 0^{\omega} \ne z 0^{\omega}$ after the input 
endmarker $\#$ was read. 
 \ \ \  \ \ \ $\Box$

\bigskip

\noindent {\bf Remark.} Although the elements of $V$ act injectively on 
$\{0,1\}^*$, the dpda in the proof of Lemma \ref{LEMdcflL_Z} is not 
injective (according to the definition of injective, or ``reversible'',
dpda in \cite{injdpda}). For example, the two configurations 
$[q_0, r \bot]$ and $[q_0, r 0^m \bot]\,$ (with $r 0^m \in {\rm domC}(b)$,
$m \ge 1$) on input $b$, yield the same next configuration
$\,b\,[q_0,\,b(r 0^m)\,\bot]$.

\medskip

\begin{pro}\!\!\!.  \label{PROwpVcyc4DCF}

\smallskip

\noindent {\small \rm (1)} \ There exists a finite monoid generating set 
$\,\Gamma_{\! H}$ of $V$ for which $ \ {\sf cowp}_{\Gamma_{\!H}\!}(V)$ $=$
$({\sf cowp}_{\Gamma_{\!H}\!}(V))^{\rm rev}$, and 

\smallskip

\hspace{1,0cm} ${\sf cowp}_{\Gamma_{\!H}\!}(V) \ \in \ $ 
${\sf cyc}(\cup_{4\,} {\sf DCFL})$ $\,\cap\,$
${\sf cyc}(\cup_{4\,} {\sf DCFL}^{\rm rev})$.

\medskip

\noindent {\small \rm (2)} \ For any finite subset $\,B = B^{-1}$
$\subseteq V$, generating a subgroup $\,\langle B \rangle\,$ 
$\subseteq\,V$, we have:

\smallskip
 
\hspace{1,0cm} ${\sf cowp}_B(\langle B \rangle)$  $\,\in\,$
$\,{\sf cyc}(\cup_{\ell\,} {\sf DCFL}^{\rm rev})$, 

\smallskip

\noindent for some constant $\,\ell \,\le\, 2^{{\rm maxlen}(B)}$.
\end{pro} 
{\sc Proof}.
(1) We use the {\em Higman generators} of $V$ as in Cor.\ \ref{CORwpVrev}.  
So, ${\rm maxlen}(\Gamma_{\! H}) = 2$. 
Every element of $\Gamma_{\! H}$ is an involution, so 
$\,\Gamma_{\! H} = \Gamma_{\! H}^{\ \pm 1}$ is a monoid generating set. 
The empty word of generators never belongs to the co-word-problem, so we
consider only words $w \in \Gamma_{\! H}^{ \ +}$.  

\smallskip

\noindent By Lemma \ref{LEMcompdepthV}(4): 
$\,w \in {\sf cowp}_{\Gamma_{\! H}\!}(V)$ \ iff \ there exists $k \in$
$[0,n]$ such that $|z_k| \le 2$ $=$ ${\rm maxlen}(\Gamma_{\! H})$, 
and $\,a_k \ldots a_1 a_n \ldots a_{k+1}(z_k 0^{\omega})\,\ne\,$
$z_k 0^{\omega}$. (If $k = 0$ then $\,a_k \ldots a_1 = \e$; if $k = n$
then $a_n \ldots a_{k+1} = \e$.)
In case $\,|z_k| = 1\,$ we replace $z_k$ by $z_k 0$ so that $|z_k| = 2$;
indeed, $z_k 0 0^{\omega} = z_k 0^{\omega}$.
The case $\,|z_k| = 0\,$ does not occur because $\e$ does not appear in the 
table of any generator of $V$.  \ Therefore,

\smallskip

\hspace{0,5cm} $w \in {\sf cowp}_{\Gamma_{\! H}\!}(V)$ \ \ iff
 \ \ $(\exists k \in [0,n])(\exists z \in \{00, 01, 10, 11\})[
\,a_k \ldots a_1 a_n \ldots a_{k+1}(z 0^{\omega})\,\ne\, z 0^{\omega}\,]$

\smallskip

\hspace{0,5cm} iff \ \ there exists $\,u \in {\sf cyc}(w)\,$ such that
$\,(\exists z \in \{00, 01, 10, 11\})[\,u \in L_z\,]$,

\smallskip

\noindent  where $ \ L_z = \{u \in \Gamma_{\! H}^{ \ *} : $ 
$u(z 0^{\omega}) \ne z 0^{\omega}\}\,$ as in Lemma \ref{LEMdcflL_Z}. 
 \ Therefore: 

\smallskip

\hspace{0,5cm} ${\sf cowp}_{\Gamma_{\! H}\!}(V)$ \ $\subseteq$ 
 \ ${\sf cyc}(L_{00})$ $\cup$  ${\sf cyc}(L_{01})$  $\cup$  
     ${\sf cyc}(L_{10})$ $\cup$  ${\sf cyc}(L_{11})$. 

\smallskip

\noindent Moreover, $\,L_z \,\subseteq\, {\sf cowp}_{\Gamma_{\! H}\!}(V)$.
And by Lemma \ref{LEMcycClos}(2), $\cup(.)$ and ${\sf cyc}(.)$ commute, so

\smallskip 

\hspace{0,5cm} ${\sf cyc}(L_{00})$ $\cup$  ${\sf cyc}(L_{01})$  $\cup$
     ${\sf cyc}(L_{10})$ $\cup$  ${\sf cyc}(L_{11})$
 $ \ = \ $  ${\sf cyc}(L_{00} \cup L_{01} \cup L_{10} \cup L_{11})$.

\smallskip

\noindent Therefore,

\smallskip

\hspace{0,5cm} $L_{00} \cup L_{01} \cup L_{10} \cup L_{11}$ \ $\subseteq$
 \ ${\sf cowp}_{\Gamma_{\! H}\!}(V)$ \ $\subseteq$
 \ ${\sf cyc}(L_{00} \cup L_{01} \cup L_{10} \cup L_{11})$.

\smallskip 

\noindent Since the word problem and the co-word-problem of a group are 
cyclically closed we finally obtain: 

\smallskip 

\hspace{0,5cm} ${\sf cowp}_{\Gamma_{\! H}\!}(V)$  $\,=\,$ 
${\sf cyc}(L_{00} \cup L_{01} \cup L_{10} \cup L_{11})$.

\smallskip

\noindent By Lemma \ref{LEMdcflL_Z} the sets $L_z$ are in 
${\sf DCFL}^{\rm rev}$, hence 
$L_{00} \cup L_{01} \cup L_{10} \cup L_{11}$ $\in$ 
$\,\cup_{4\,}{\sf DCFL}^{\rm rev}$. So 
$\,{\sf cowp}_{\Gamma_{\! H}\!}(V)$ $\in$ 
$\,{\sf cyc}(\cup_{4\,}{\sf DCFL}^{\rm rev})$.

\smallskip

By Cor.\ \ref{CORwpVrev}, ${\sf cowp}_{\Gamma_{\! H}\!}(V)$ is closed
under reversal; and over $\Gamma_{\! H}$, $\,L_z = L_z^{\rm \,rev}$.
And ${\sf cyc}(.)$, $(.)^{\rm rev}$, and $\cup$ commute.  
Hence, $\,{\sf cowp}_{\Gamma_{\! H}\!}(V)$ $\in$ 
$\,{\sf cyc}(\cup_{4\,}{\sf DCFL})$ \ (without reversal).  

\medskip

\noindent (2) In a similar way as for $\Gamma_{\!H}$ one shows that for any 
finite set $B = B^{-1} \subseteq V$: 

\smallskip

\hspace{0,5cm} ${\sf cowp}_B(\langle B \rangle)$ $ \ = \ $ 
${\sf cyc}\big(\bigcup \{L_z : z \in \{0,1\}^{{\rm maxlen}(B)}\}\big)$ 

\medskip

\hspace{0,5cm} $\in\,$ ${\sf cyc}(\cup_{\ell\,} {\sf DCFL}^{\rm rev})\,$, 
 
\medskip

\noindent where $\,\ell \,\le\, 2^{{\rm maxlen}(B)}$.
 \ \ \  \ \ \ $\Box$

\bigskip

\noindent {\bf Remarks.}

\smallskip

\noindent (1) \ For Prop.\ \ref{PROwpVcyc4DCF}, it does not matter in what 
order the operations $\cup$, $(.)^{\rm rev}$, and ${\sf cyc}(.)$ are 
written in $\,{\sf cyc}(\,\cup_{\ell\,} {\sf DCFL}^{\rm rev})$.
Indeed, these operations commute two-by-two.

\smallskip

\noindent (2) \ Because of Prop.\ \ref{DCFLcyc}, the operation ${\sf cyc}(.)$
cannot simply be removed in Prop.\ \ref{PROwpVcyc4DCF}.

\noindent Since ${\sf cyc}(.)$ and ${\sf co}(.)$ do not commute, by Prop.\
\ref{LEMcycClos}(1.b), ${\sf co}(.)$ cannot simply be removed either.

\medskip

\begin{cor} \label{CORwpVn2}\!\!.

\smallskip

\noindent {\small \rm (1)} The word problem of $V$ over any finite 
generating set is in $\,{\sf DTime}(n^2)$.

\smallskip

\noindent {\small \rm (2)} The word problems of the {\em Higman-Thompson 
groups} $G_{n,r}$, $T_{n,r}$, and $F_{n,r}$ (in particular $F$), over any 
finite generating set are in ${\sf DTime}(n^2)$.
\end{cor}
{\sc Proof.} (1) We mentioned that $\,\cup_{\ell\,} {\sf DCFL}\,$ and 
$\,\cup_{\ell\,} {\sf DCFL}^{\rm rev}\,$ are in ${\sf DTime}(n)$. And if
$L \in {\sf DTime}(n^d)$ then ${\sf cyc}(L) \in {\sf DTime}(n^{d+1})$; 
this follows from the fact that a word $w$ can be cyclicly permuted in at 
most $|w|$ different ways. Since ${\sf DTime}(n^2)$ is closed under 
complementation, it contains a word problem iff it contains the 
co-word-problem.

(2) This follows from the fact that $G_{n,r}$, etc., are finitely generated 
subgroups of $V\,$ \cite{Hig74}.
 \ \ \ $\Box$

\subsection{The class {\sf logDCFL} }

By Prop.\ \ref{PROwpVcyc4DCF} the word problem of $V$ and of all of its 
finitely generated subgroups belongs to the closure of {\sf DCFL} under the
operations ${\sf co}(,)$, $(.)^{\rm rev}$, $\cup$, and ${\sf cyc}(.)$. 
This class is contained in {\sf logDCFL}, which is also closed under these
operations, as we will see in the next Proposition.  

By definition, {\sf logDCFL} is the set of languages that can
be reduced to a language in {\sf DCFL} by a many-one log-space reduction. 
I.e., for $L \subseteq A^*$ we have $L \in {\sf logDCFL}$ \ iff \ there 
exists $K \subseteq B^*$ with $K \in {\sf DCFL}$, and there exists a 
total function $f\!: A^* \to B^*\,$ (the {\em reduction function}) that is 
computable by a deterministic log-space multitape Turing machine, such that
$L = f^{-1}(K)$; see \cite{Johnson}.

Similarly, ${\sf log}_{(1:1)}{\sf DCFL}$ denotes the closure of {\sf DCFL} 
under one-one log-space reduction, i.e., one uses {\em injective} log-space 
computable reduction functions; see \cite{Johnson}. 

\smallskip

\begin{pro} \label{PROPlogDCFLclos} 
 \ The classes $\,{\sf logDCFL}$ and $\,{\sf log}_{(1:1)}{\sf DCFL}$ are 
closed under the operations of reversal, complementation, finite 
intersection, finite union, and cyclic permutation.
\end{pro}
{\sc Proof.}
(1) \ {\sf Claim.} If $L \in {\sf logDCFL}$ then $L^{\rm rev} \in$
${\sf logDCFL}$. The same holds for ${\sf log}_{(1:1)}{\sf DCFL}$.

\smallskip

\noindent {\sf Proof of the Claim.} We define the reduction function $f$
of $L^{\rm rev}$ to $L$ by $f(x)$ $=$ $x^{\rm rev}$.
Then indeed $y \in L^{\rm rev}$ iff $f(y) = y^{\rm rev}$ $\in$ $L$. 
The function $f$ can be computed in space complexity 0, and is injective.

\medskip

\noindent (2) \ {\sf Claim.} If $L_1, \,\ldots\, , L_k \in {\sf logDCFL}$ 
for some $k \ge 2$, then $\,L_1 \cap \,\ldots\, \cap L_k \in$ 
{\sf logDCFL}. The same holds for ${\sf log}_{(1:1)}{\sf DCFL}$. 

\smallskip

\noindent {\sf Proof of the Claim.} The claim follows by induction from 
the case where $k = 2$.  
For $i = 1, 2$, let $L_i \in {\sf logDCFL}$, $L_i \subseteq A^*$, and
$L_i = f_i^{-1}(L_i^{(o)})$ $\subseteq$ $A^*$, where 
$f_i\!: A^* \to B^*$ is log-space computable and 
$L_i^{(o)} \in {\sf DCFL}$, $\,L_i^{(o)} \subseteq B^*$.  For a new letter 
$\# \not\in B$, let $L^{(o)} =$ $L_1^{(o)} \# L_2^{(o)}\,$ (the 
{\em marked concatenation}), which belongs  to {\sf DCFL}; see 
\cite[{\small Sect.\ 11.3 Ex.\ 2}]{Harrison}.

We reduce $L_1 \cap L_2\,$ to $L^{(o)}$ by the function $f\!: x \in A^*$ 
$\longmapsto$ $f_1(x)\,\#\,f_2(x) \,\in\, B^*\#B^*$. 
Then $\,x \in L_1 \cap L_2$ \ iff \ $f_1(x) \in L_1^{(o)}$ and 
$f_2(x) \in L_2^{(o)}$ \ iff \ $f_1(x)\,\#\,f_2(x) \in$
$L_1^{(o)} \# L_2^{(o)}$.  
It is straightforward to see that $f$ is computable in log-space.

It is easy to prove that if $f_1$ or $f_2$ is injective then so is $f$.

\medskip

\noindent (3) \ {\sf Claim.} If $L \in {\sf logDCFL}$ then 
${\sf co}(L) \in {\sf logDCFL}$. The same holds for
${\sf log}_{(1:1)}{\sf DCFL}$.

\smallskip

\noindent {\sf Proof of the Claim.} Let $L = f_o^{-1}(L_o) \subseteq A^*$ 
where $L_o$ $\subseteq B^*$, $\,L_o \in {\sf DCFL}$, and 
$f_o\!: A^* \to B^*$ is log-space computable.  
Then $\,{\sf co}(L_o) = B^* \minus L_o \in$ ${\sf DCFL}$, since 
{\sf DCFL} is closed under complementation.  
Then the same function $f_o$ reduces $A^* \minus L\,$ to 
$\,f_o^{-1}(B^* \minus L_o)$, since $\,f_o^{-1}(B^* \minus L_o)$ $\,=\,$ 
$A^* \minus f_o^{-1}(L_o)$. 

\medskip

\noindent (4) \ {\sf Claim.} \ {\sf logDCFL} and 
${\sf log}_{(1:1)}{\sf DCFL}$  are  closed under finite union.

\smallskip

\noindent {\sf Proof of the Claim.} This follows from (2) and (3).  

\medskip

\noindent (5) \ {\sf Claim.} If $L \in {\sf logDCFL}$ then ${\sf cyc}(L)$
$\in {\sf logDCFL}$. The same holds for ${\sf log}_{(1:1)}{\sf DCFL}$.

\smallskip

\noindent {\sf Proof of the Claim.} Let $L = f_0^{-1}(L_0) \subseteq A^*$ 
where $L_0 \subseteq B^*$, $\,L_0 \in {\sf DCFL}$, and $f_0\!: A^* \to B^*$
is log-space computable. Let $\# \not\in B$ be a new letter. We use the fact
that {\sf DCLF} is closed under {\em marked Kleene star};  for a language 
$K \subseteq B^*$, the marked Kleene star is $(K\#)^*\,$; see 
\cite[{\small Sect.\ 11.3 Ex.\ 2}]{Harrison}.  

Let $L_c = (B^*\# \minus L_0 \#)^*$, which  belongs to {\sf DCFL}; indeed,
$B^* \minus L_0 \in {\sf DCFL}$ by closure under complementation, and 
$\,(B^*\# \minus L_0 \#)^* = ((B^* \minus L_0)\#)^* \in {\sf DCFL}$ by
closure under marked Kleene star. We will reduce $\,{\sf co}({\sf cyc}(L))$
to $\,L_c$.

Let $\kappa$ be the {\em one-step cyclic permutation}; i.e., for any 
$x_0 x_1 \,\ldots\, x_{n-1} \in A^*$ with $x_i \in A$ for 
$0 \le i$ $\le$ $n\!-\!1$, we define $\,\kappa(x_0x_1 \,\ldots\, x_{n-1})$ 
$=$ $x_1 \,\ldots\, x_{n-1} x_0$.  The reduction function $f$: $A^* \to$ 
$(B \cup \{\#\})^*\,$ is defined by 

\smallskip

\hspace{0,5cm} $f(x) \,=\, $
$f_0(x)\,\#\,f_0(\kappa(x))\,\#\,f_0(\kappa^2(x))\,\#\,$ $ \,\ldots \,$  
$f_0(\kappa^j(x))\,\#\,$  $ \,\ldots \,$ $f_0(\kappa^{n-1}(x))\,\#\,$
 \ \ $\in$ \ $(B^*\#)^*$,

\smallskip

\noindent where $n = |x|$. Then $f$ reduces $\,{\sf co}({\sf cyc}(L))$ $=$
$A^* \minus {\sf cyc}(L)\,$ to $\,L_c =$ $(B^*\# \minus L_0\#)^*$. \\  
Indeed, $\,x \in A^* \minus {\sf cyc}(L)$ \ iff 
 \ $(\forall j \in [0, n\!-\!1])[\,\kappa^j(x) \not\in L\,]$ 
 \ iff \    
$(\forall j \in [0,n\!-\!1])[\,f_0(\kappa^j(x)) \in B^* \minus L_0\,]$ \\  
iff \ \  
$f_0(x)\,\#\,f_0(\kappa(x))\,\#\,f_0(\kappa^2(x))\,\#\,$ $\,\ldots\,$
$f_0(\kappa^j(x))\,\#\,$  $\,\ldots\,$ $f_0(\kappa^{n-1}(x))\,\#\,$
 \ $=$ \ $f(x)$ \ $\in$ \ $(B^*\# \minus L_0\#)^*$.

\smallskip

\noindent To show that $f$ is computable in log-space we express $f$ as the 
composite of two functions:

\medskip

\hspace{0,5cm} $K: \ x \in A^* \,\longmapsto\, $
$x\,\#\,\kappa(x)\,\#\,\kappa^2(x)\,\#\,$ $ \,\ldots \,$
$\kappa^j(x)\,\#\,$  $ \,\ldots \,$ $\kappa^{|x|-1}(x)\,\#\,$
 \ \ $\in \ (A^*\#)^*\,$;

\bigskip

\hspace{0,5cm} $F_{\!o}: \ z^{(0)}\,\#\,z^{(1)}\,\#\,z^{(2)}\,\#\,$ 
$\,\ldots\,$ $z^{(j)}\,\#\,$  $\,\ldots\,$ $z^{(m-1)}\,\#\,$ 
 \ \ $\in \ (B^*\#)^*$ \ \ \ (for any $m \ge 0$) 

\smallskip

\hspace{1,5cm} $\longmapsto$
 \ $f_0(z^{(0)})\,\#\,f_0(z^{(1)})\,\#\,f_0(z^{(2)})\,\#\,$ $ \,\ldots \,$
$f_0(z^{(j)})\,\#\,$  $ \,\ldots \,$ $f_0(z^{(m-1)})\,\#\,$
 \ \ $\in \ (B^*\#)^*$. 

\medskip

\noindent So $f(.) = F_{\!o} \circ K(.)$.
Since the composite of two log-space computable functions is 
log-space computable (see e.g.\ \cite[{\small Lemma 13.3}]{HU}), it
suffices to show that $K(.)$ and $F_{\!o}(.)$ are log-space computable.
For $F_{\!o}(.)$ this is easy. 
The function $K$ is obviously injective, and one proves easily that if 
$f_0$ is injective then so is $F_o$.

To compute $K(.)$ in log-space we consider a Turing machine that 
successively computes the words $\kappa^j(x)\,\#$\, on input $x$ for 
$j = 0, 1, 2,$ $\,\ldots\,, n\!-\!1\,$ (where $n = |x|$). 
The number $j$ $\in [0, n\!-\!1]$, written in binary, is stored in space 
$\log n$, initialized to $0$.
To compute $\,x_0x_1 \,\ldots\, x_{n-1}$ $\longmapsto$
$x_j \,\ldots\, x_{n-1} x_0 \,\ldots\, x_{j-1}\#$, a second log-space counter
$i$ is used, written in binary and initialized to $j$. The machine starts at 
the left end of $x$ and moves right while decrementing $i$. When $i=0$, the 
head is on $x_j$.  Now $x_j \,\ldots\, x_{n-1}$ is copied to the output, 
until the input head reaches the right endmarker of the input tape. 
Next, the input head is moved left to the left endmarker of the input tape, 
and $i$ is re-initialized to $j$. 
The head then moves right and $x_0 \,\ldots\, x_{j-1}$ is copied to the
output, while $i$ is being decremented. When $i=0$, $\#$ is printed on the
output; $\kappa^j(x)\,\#$ has been produced. Now $j$ is incremented to 
$j\!+\!1$, the input head moves back to the left end, and the process 
repeats with the new value of $j$, unless $j = n$. If $j = n$, the 
computation of $K(x)$ is complete, and stops. 

In summary, we have found a log-space reduction of 
$\,{\sf co}({\sf cyc}(L))\,$ to $\,L_c =$ $(B^*\# \minus L_0\#)^*$ $\in$ 
${\sf DCFL}$. Hence $\,{\sf co}({\sf cyc}(L))$ $\in$ {\sf logDCFL}. 
By closure under complementation, proved in (3), we obtain 
$\,{\sf cyc}(L)$ $\in$ {\sf logDCFL}. 
We saw that the same applies to ${\sf log}_{(1:1)}{\sf DCFL}$.
 \ This proves the Claim.
 \ \ \  \ \ \ $\Box$

\bigskip

\noindent {\sf Notation:} If ${\cal C}$ is a set of languages then
$\,\cap_{\ell\,} {\cal C}\,$ denotes the set of languages that are the 
intersection of $\le \ell$ languages in $\cal C$.

\begin{pro} \label{PROPlogDCFL}
 \ \ ${\sf co}({\sf cyc}(\cup_{\ell\,} {\sf DCFL}^{\rm rev}))$  $ \ = \ $
$\cap_{\ell\,} {\sf co}({\sf cyc}({\sf DCFL}^{\rm rev}))$
$ \ \subseteq \ {\sf log}_{(1:1)}{\sf DCFL}$.
\end{pro}
{\sc Proof.} The equality follows from the commutativity relations 
between $\cup_{\ell}(.)$, ${\sf cyc}(.)$, and ${\sf co}(.)$, where 
${\sf co}(\cup_{\ell}(.))$ $=$ $\cap_{\ell}({\sf co}(.))$.
The inclusion then follows from the closure properties of {\sf logDCFL}
and ${\sf log}_{(1:1)}{\sf DCFL}$ in Prop.\ \ref{PROPlogDCFLclos}.
 \ \ \  \ \ \ $\Box$

\begin{cor}
 \ The word problem of $V$ over any finite monoid generating set, and the
word problem of any finitely generated subgroup, is in 
${\sf log}_{(1:1)}{\sf DCFL}$ \ $(\,\subseteq$ {\sf logDCFL}). 
\end{cor} 
{\sc Proof.} This follows immediately from Prop.\ \ref{PROPlogDCFL}, Cor.\
\ref{CORwpVn2}, and Prop.\ \ref{PROPlogDCFLclos}.
 \ \ \  \ \ \ $\Box$

\subsection{Final remarks and questions }

\noindent {\sf (1) The Thompson group $F$:}
 \ Shpilrain and Ushakov \cite{ShpUsh} show that the word problem of $F$ 
can be decided by a $O(n\,\log n)$-step program, 
based on the infinite presentation 
$\,\langle \{x_i : i \in \omega\}$ $:$ $\{ x_i^{-1} x_k x_i = x_{k+1} :\,$ 
$i,k \in \omega, \ {\rm and} \ i < k\} \rangle\,$ of $F$; \ here, $n$ is 
the length of the input over the infinite alphabet 
$\{x_i, x_i^{-1} : i \in \omega\}$.
This is a useful result, but it does not directly yield the time-complexity 
of the word problem of $F$ over a finite generating set, on a multitape 
Turing machine; nor does it directly give the Dehn function of $F$, which 
is quadratic \cite{Guba}.  

\medskip

\noindent {\sf (2) Lehnert's conjecture} \cite{Lehnert,BMN}
says that every finitely generated group with word problem in {\sf coCFL} is
isomorphic to a subgroup of $V$.\footnote{ \ Lehnert \cite{Lehnert} stated 
a somewhat different conjecture, but Bleak, Matucci and Neunh\"offer 
\cite{BMN} proved that this is equivalent to the present form of the 
conjecture.}
 \ By Prop.\ \ref{PROwpVcyc4DCF}(2) this would imply that every finitely 
generated group with {\sf coCFL} word problem actually has its word 
problem in $\,{\sf co}({\sf cyc}(\cup {\sf DCFL}^{\rm rev}))$ ($\subseteq$
${\sf DTime}(n^2)$).  

\medskip

\noindent {\sf (3)} We did {\sl not} prove that for every two finitely 
generated groups $G_1 = \langle B_1 \rangle$ and $G_2 = \langle B_2 \rangle$
in general, if $G_1 \subseteq G_2$ and ${\sf wp}_{B_2}(G_2)$ $\in$  
${\sf co}({\sf cyc}(\cup {\sf DCFL}^{\rm rev}))$, then 
${\sf wp}_{B_1}(G_1)$ $\in$ 
${\sf co}({\sf cyc}(\cup {\sf DCFL}^{\rm rev}))$.

\medskip

\noindent {\sf (4)} It remains open whether 
${\sf cyc}(\cup {\sf DCFL}^{\rm rev})$ is a strict 
subset of ${\sf cyc}({\sf CFL})\,$ $(\subseteq {\sf CFL})$. We know that 
$\cup {\sf DCFL}^{\rm rev}$ $\subsetneqq$ ${\sf CFL}$, and 
${\sf cyc}(\cup {\sf DCFL}^{\rm rev})$ $\subseteq$ 
${\sf DTime}(n^2)\,$ (whereas ${\sf CFL}$ $\subseteq$ ${\sf DTime}(n^2)$ 
would be surprising). 

\medskip

\noindent {\sf (5)} A few more open questions: 

\smallskip

\noindent $\,\circ$ Is the word problem of $V$ (or of $F$) {\em complete} in 
${\sf log}_{(1:1)}{\sf DCFL}$ for one-one log-space reduction  
\cite{Sudborough,Lohrey}?  

\smallskip

\noindent $\,\circ$ Does ${\sf wp}(V)$ reduce to ${\sf wp}(F)\,$?
(True if ${\sf wp}(F)$ were complete in ${\sf log}_{(1:1)}{\sf DCFL}$.)

\smallskip

\noindent $\,\circ$ Does ${\sf cowp}(V)$ reduce to ${\sf wp}(V)\,$? (True
if ${\sf wp}(V)$ were complete in ${\sf log}_{(1:1)}{\sf DCFL}$ .) 

\medskip

\section{Appendix: Push-down automata}

The material of this Appendix goes back to the early 1960s, and is exposited
in many more books and articles than the ones cited here.

A push-down automaton (pda) is a structure ${\cal A}$ $=$
$(Q, A,$  $\Sigma,$ ${\cal T},$ $q_0,$  ${\rm s}_0,$  $Q_{\rm a})$, where
$Q$ (state set), $A$ (input alphabet), $\Sigma$ (stack alphabet), ${\cal T}$
(set of transitions), and $Q_{\rm a} \subseteq Q$ (set of accept states), 
are finite sets; $q_0 \in Q$ is the start state; and ${\rm s}_0 \in \Sigma^+$
is the initial content of the stack.

The current {\em configuration} of ${\cal A}$ (a.k.a.\ the instantaneous
description) is of the form $ \ w \ [q, s]\,$, where $q \in Q$ is the current
state,  $s \in \Sigma^*$ is the current stack content, and $w \in A^*$ is 
the input that has been read so far. 
A pda has one {\em start configuration}, namely $[q_0,s_0]$, where $q_0$ and
$s_0$ are as above. No input has been read at this point, so $w$ is $\e$ in
the start configuration.
An {\em accept configuration} is of the form $\,w\, [q, s]$, such that 
$q \in Q_{\rm a}$. (See the $\exists$- and the $\forall$-acceptance rules
below.)

\smallskip

\noindent Remark: Our definition of configuration is different from the one
in the literature \cite{HU, Harrison, Ginsburg}. In these books a 
configuration is of the form $(q, x, s) \in$ $Q \x A^* \x \Sigma^*$, where
$q$ and $s$ are the same as for us, but $x$ is a {\em future} input. 
For us, $w$ in $\,w\,[q,s]\,$ is the past input that has been read.  

\medskip

A {\em transition} in ${\cal T}$ has the form  
$ \ (q,s) \stackrel{a}{\to} (p,s)$, where 
$(q, s) \in Q \x \Sigma^+$, $(p,s') \in Q \x \Sigma^*$, and $a \in$ 
$A \cup \{\e\}$. When $a = \e$, this is called an $\e$-transition: the state 
and the stack may change, but no next input letter is being read (either the
next input letter is not yet there, or it is there but is not yet being read).
There is no transition on an empty stack; in a transition as above, 
$s \in \Sigma^+$. A pda has a finite set ${\cal T}$ of transitions
(called the {\em transition table}). 

The transition $(q,s) \stackrel{a}{\to} (p,s')$ is {\em applicable} to a 
configuration $w\,[r',t']$ $\in$  $Q \x \Sigma^+\,$ iff $\,r' = q$, $\,s$ 
is a prefix of $t'$, and either $a \in A$ and the next input letter is $a$,
or $a = \e$ (and then there is no requirement on the input). When this 
transition is {\em applied to the configuration} $w\,[q, s t]$, the next 
configuration is $wa \,[p,s't]\,$ (where $a$ can be $\e$, in an 
$\e$-transition).
We extend the transition notation to configurations: when the transition
$(q, s) \stackrel{a}{\to} (p,s')\,$ is applied to the configuration
$\,w \,[q, s t]$, we write
$\,w \,[q, s t] \stackrel{a}{\to} \,wa\,[p, s' t]$.

Our transitions are a little more general than the ones commonly used in 
the literature, but they do not lead to the acceptance of more languages 
\cite{HU, Harrison, Ginsburg, Sipser}. 

A {\em computation} of a pda ${\cal A}$ on input $w =$ 
$a_1 a_2 \,\ldots\, a_n \in A^*$ is a sequence of configurations and 
applications of transitions

\medskip

\hspace{0.5cm}  $[q,s_0]$
$\,\stackrel{\e^*}{\longrightarrow}\, [q_1', s_1']$
$\,\stackrel{a_1}{\to}\, a_1\,[q_1, s_1]$
$\,\stackrel{\e^*}{\longrightarrow}\, a_1\,[q_2',s_2']$
$\,\stackrel{a_2}{\to}\, a_1 a_2\,[q_2, s_2]$
$\,\stackrel{\e^*}{\longrightarrow}  \ \ \ldots \ \ $

\smallskip

\hspace{0.7cm}  $ \ \ \ldots \ \ \ $
$\stackrel{\e^*}{\longrightarrow}\, a_1 a_2 \ldots a_{n-1}\,[q_n',s_n']$
$\,\stackrel{a_n}{\to}\, a_1 a_2 \ldots a_{n-1}a_n\,[q_n,s_n]$
$\,\stackrel{\e^*}{\longrightarrow}\, a_1 a_2 \ldots a_{n-1}a_n\,[p,s]$,

\medskip

\noindent where $\,\stackrel{\e^*}{\longrightarrow}\,$ denotes a (possibly
empty) sequence of $\e$-transitions.

\medskip

The language accepted by a pda depends on the acceptance rule.

The $\exists$-acceptance rule for a pda ${\cal A'}$ is as follows:
 \ $w \in A^*$ is {\em accepted} iff {\em there exists} a computation of 
${\cal A'}$ on input $w$ that reads the whole input $w$, starting with the 
start configuration $[q_0,s_0]$, and ends with a configuration of the form 
$\,w\,[q,t]$ for some $q \in Q_{\rm a}$, $t \in \Sigma^*$ (depending on $w$). 
A pda ${\cal A'}$ that uses the $\exists$-acceptance rule is called a 
$\exists$pda. 

The $\forall$-acceptance rule for a pda ${\cal A}$ is as follows:
 \ $w$ is {\em accepted} iff {\em every} computation of ${\cal A}$ on input
$w$ reads the whole input $w$, starting with the start configuration 
$[q_0,s_0]$, and ends in a configuration $\,w\,[q,t]$ for some 
$q \in Q_{\rm a}$, $t \in \Sigma^*$ (depending on $w$).
A pda ${\cal A}$ that uses the $\forall$-acceptance rule is called a 
$\forall$pda.

The same pda can be used as a $\exists$pda or a $\forall$pda; it depends on
the acceptance rule.

\smallskip

A {\em deterministic pda} (dpda) is a pda such that in every configuration, 
at most one transition is applicable. The $\forall$-acceptance rule is then 
the same as the $\exists$-acceptance rule.

\medskip

We call a stack symbol $\bot \in \Sigma$ a {\em bottom marker} \ iff
{\small \rm (1)} \ the start configuration is $[q_0, \,p_0 \bot]$
for some $p_0 \in (\Sigma \minus \{\bot\})^*$;
 \ {\small \rm (2)} \ every transition with left-side $(q, s\bot)$
has a right-side $(p, s'\bot)$, for some $q,p \in Q$ and
$s, s' \in (\Sigma \minus \{\bot\})^*$;
 \ {\small \rm (3)} \ every transition with left-side $(q, s)$
where $s \in (\Sigma \minus \{\bot\})^+$, has a right-side $(p,s')$ for 
some $s' \in (\Sigma \minus \{\bot\})^*$. \ (Conditions {\small (2)} and
{\small (3)} say that $\bot$ is never erased nor written.)

\medskip

It is well known that the $\exists$pda$s$ accept exactly the languages in 
{\sf CFL}, and the $\forall$pda$s$ accept exactly the languages in
{\sf coCFL}. By definition, {\sf DCFL} is the set of languages accepted
by dpda$s$. We have $\,{\sf DCFL} \subsetneqq {\sf CFL} \cap {\sf coCFL}$; 
for ``$\ne$'', the set of palindromes and the set of non-palindromes each 
has a fairly simple CF grammar. 


\bigskip

\noindent {\bf Acknowledgement.} I would like to thank the referee for 
many useful comments. 




\medskip

\noindent {\small birget@camden.rutgers.edu}

\end{document}